\documentclass{amsart}

\usepackage{amssymb,amsmath}

\usepackage{color}
\theoremstyle{plain}
\newtheorem{theorem}{Theorem}[section]
\newtheorem*{theorem*}{Theorem}

\theoremstyle{remark}

\numberwithin{equation}{section}
\theoremstyle{definition}
\newtheorem{definition}[theorem]{Definition}
\numberwithin{equation}{section}
\numberwithin{equation}{section}

\newcommand\RR{\mathbb{R}}
\newcommand\dir{{\rm{\,d}}\rho}
\newcommand\dil{{\rm{\,d}}\lambda}
\newcommand\di{{\rm{\,d}}}
\newcommand\nep{{\rm{e}}}
\newcommand\CZ{Calder\'on--Zygmund }

\begin{document}

\title[Hardy spaces and Riesz transforms]{Hardy spaces and Riesz transforms \\on a Lie group of exponential growth }



\thanks{This work is partially supported by the project "Harmonic analysis on continuous and discrete structures" funded by Compagnia di San Paolo
		(Cup E13C21000270007). The second-named author is a member of the
Gruppo Nazionale per l'Analisi Matematica, la Probabilit\`a e le loro Applicazioni
(GNAMPA) of the Istituto Nazionale di Alta Matematica (INdAM). The first-named author is grateful to the Politecnico di Torino for its hospitality during several visits which made this work possible.}

\author[P. Sj\"ogren and M. Vallarino]
{Peter Sj\"ogren and Maria Vallarino}

\address{Peter Sj\"ogren:
Mathematical Sciences,
University of Gothenburg and  Mathematical Sciences,
Chalmers\\
S-412 96 G\"oteborg\\Sweden
}
\email{peters@math.chalmers.se}

\address{Maria Vallarino:
Dipartimento di Scienze Matematiche "Giuseppe Luigi Lagrange" -
\\ Politecnico di Torino\\
Corso Duca degli Abruzzi, 24 \\ 10129 Torino \\ Italy} \email{maria.vallarino@polito.it}

\begin{abstract}
Let $G$ be the Lie group ${\Bbb{R}}^2\rtimes {\Bbb{R}}^+$ endowed with the Riemannian symmetric space structure. Take  a distinguished basis  $X_0,\, X_1,\,X_2$ of left-invariant vector fields of the Lie algebra of $G$, and consider the Laplacian $\Delta=-\sum_{i=0}^2X_i^2$ and the first-order Riesz transforms $\mathcal R_i=X_i\Delta^{-1/2}$, \hskip3pt $i=0,1,2$. We first show that the atomic Hardy space $H^1$ in $G$ introduced by the authors in a previous paper does not admit a characterization in terms of the Riesz transforms $\mathcal R_i$. It is also proved that two of these Riesz transforms are bounded from $H^1$ to $H^1$.

\end{abstract}

\maketitle

\section{Introduction}\label{intro}
Let $G$ be the Lie group $\RR^2\rtimes \RR^+$ where the product rule is the following:
$$(x_1,x_2,a)\cdot(x'_1,x'_2,a')=(x_1+a\,x'_1,x_2+a\,x'_2,a\,a') $$
for $(x_1,x_2,a),\,(x'_1,x'_2,a')\in G$. We shall denote by $x$ the point $(x_1,x_2,a)$,
and it will be convenient to write
\begin{equation*}
   |x|=\sqrt{x_1^2+x_2^2} \qquad \text{for} \quad \; x = (x_1,x_2,a).
\end{equation*}
The group $G$ is not unimodular; a right and a left Haar measure of $G$ are given by
$$\dir(x)=a^{-1}\,\di x_1\di x_2\di a\qquad{\rm{and}}\qquad \dil(x)=a^{-3}\,\di x_1\di x_2\di a\,,$$
respectively. The modular function is thus $\delta(x)=a^{-2}$. Throughout this paper, unless explicitly stated, we use the right measure $\rho$ on $G$ and denote by $L^p$, $\|\cdot\|_p$ and $\langle\cdot,\cdot \rangle$  the $L^p$-space, the $L^p$-norm and the $L^2$-scalar product with respect to  $\rho$.

The group $G$ has a Riemannian symmetric space structure, and the corresponding metric, which we denote by $d$, is that of the three-dimensional hyperbolic half-space. The metric $d$ is invariant under left translation and given by
\begin{equation}\label{metrica}
\cosh r(x)=\frac{a+a^{-1}+a^{-1}|x|^2}{2}\,,
\end{equation}
where $r(x)=d(x,e)$ denotes the distance of the point $x$ from the identity $e=(0,0,1)$ of $G$. The measure of a hyperbolic ball $B_r$, centred at the identity and of radius $r$, behaves like
$$\lambda(B_r)=\rho(B_r)\sim  \begin{cases}
r^3&{\rm{if}}~r<1\\
\nep^{2r}&{\rm{if}}~r\geq 1\,.
\end{cases}
$$

Thus $G$ is a group of {\emph{exponential volume growth}}, and so nondoubling. In this context, the classical \CZ theory and the classical definition of the atomic Hardy space (see \cite{CW1, CW2, St}) do not apply. But Hebisch and Steger \cite{HS} have constructed a \CZ theory which applies to some nondoubling spaces, in particular to our space $(G,d,\rho)$. The main idea is to replace the family of balls in the classical \CZ theory by a suitable family of parallelepipeds which we call \emph{\CZ sets}. The definition appears in \cite{HS} and implicitly in \cite{GS}, and reads as follows.

\begin{definition}\label{Czsets}
A Calder\'on-Zygmund set is a parallelepiped $P=[x_1-L/2,x_1+L/2]\times[x_2-L/2,x_2+L/2]\times [a\nep^{-r},a\nep^r]$, where $L>0$, $r>0$ and $(x_1,x_2,a)\in G$ are related by
$$\nep^{2}a\,r\leq L< \nep^{8 }a\,r\qquad{\rm{if~}}r<1\,,$$
$$a\,\nep^{2r}\leq L< a\,\nep^{8r}\qquad{\rm{if~}}r\geq 1\,.$$
The point $(x_1,x_2,a)$ is the center of $P$, and we call $r$ the parameter of $P$.
\end{definition}

We let $\mathcal P$ denote the family of all \CZ sets, and observe that  $\mathcal P$ is invariant under left translation. In \cite{HS} it is proved that every integrable function on $G$ admits a \CZ decomposition involving the family $\mathcal P$, and that a \CZ theory can be developed in this context. Using the \CZ sets, it is natural to introduce an atomic Hardy space $H^1$ on the group $G$, as follows (see \cite{V} for details).

For $1<p\leq \infty$, a $(1,p)$-atom is a function $A$ in $L^1$ such that
\begin{itemize}
\item [(i)] $A$ is supported in a \CZ set $P$;
\item [(ii)]$\|A\|_{p}\leq \rho(P)^{-1+1/p}\,;$
\item [(iii)]$\int A\dir =0$\,.
\end{itemize}
The atomic Hardy space is now defined in a standard way.
\begin{definition}
The atomic Hardy space $H^{1,p}$ is the space of all functions $f$ in $ L^1$ which can be written as $f=\sum_j \lambda_j\, A_j$, where the $A_j$ are $(1,p)$-atoms and $\lambda _j$ are complex numbers such that $\sum _j |\lambda _j|<\infty$. We denote by $\|f\|_{H^{1,p}}$ the infimum of $\sum_j|\lambda_j|$ over such decompositions.
\end{definition}
By \cite[Theorem 2.3]{V}, for any $p\in (1,\infty)$ the space $H^{1,p}$ coincides with $H^{1,\infty}$ and their norms are equivalent. In the following we shall simply denote this space by $H^1$ and its norm by $\|\cdot\|_{H^1}$.

\smallskip

The \CZ theory from \cite{HS} has turned out to be useful to study the boundedness of singular integral operators related to the distinguished Laplacian on $G$, defined as follows.

Let $X_0,\,X_1,\,X_2$ denote the left-invariant vector fields
$$X_0=a\,\partial_a,\qquad X_1=a\,\partial_{x_1},\qquad  X_2=a\,\partial_{x_2}\,,$$
which span the Lie algebra of $G$. The Laplacian $\Delta=-(X_0^2+X_1^2+X_2^2)$ is a left-invariant operator which is essentially selfadjoint on $L^2$. It is well known that the heat semigroup $\big( \nep^{-t\Delta}\big)_{t>0}$ is given by a kernel $h_t$, in the sense that $\nep^{-t\Delta}f=f\ast h_t$ for suitable functions $f$.  Let $\mathcal M_h$ denote the corresponding heat maximal operator, defined by
\begin{equation} \label{hetamaxop}
\mathcal M_h\, f(x)=\sup_{t>0} |f\ast h_t(x)|\qquad \forall x\in G\,.
\end{equation}
We then define the heat maximal Hardy space $H^{1}_{\rm{max},h}$ as the space of all functions $f$ in $ L^1$ such that $\mathcal M_h\, f$ is in $L^1$, endowed with the norm $$\|f\|_{H^{1}_{\rm{max},h}}=\|\mathcal M_hf\|_1\,.
$$
The authors proved in \cite{SV2} that $H^{1}\subset H^{1}_{\rm{max},h}$ and that this inclusion is strict.

Similarly, one can consider a maximal Hardy space $H^{1}_{\rm{max},p}$ on $G$ by means of the Poisson maximal function. It also strictly contains the atomic Hardy space, see  \cite{V2}. Thus there is no characterization of the atomic Hardy space by means of the heat or the Poisson
maximal operator, in our setting.

We shall now consider the first-order Riesz transforms associated with $\Delta$, defined by $\mathcal R_i=X_i\,\Delta^{-1/2}$, \,\,$i=0,1,2$. The associated Hardy space is
$$
H^{1}_{\rm{Riesz}}=\{f\in L^1:\,\mathcal R_if\in L^1,\,\,\,i=0,1,2\}\,,
$$
endowed with the norm
$$\|f\|_{H^{1}_{\rm{Riesz}}}=\|f\|_1+\sum_{i=0}^2\|\mathcal R_if\|_1\,.
$$
The authors proved in \cite{SV1} that the operators $\mathcal R_i$, \,\,$i=0,1,2$, are bounded from $H^1$ to $L^1$. Thus $H^{1}$ is included in $H^{1}_{\rm{Riesz}}$. We prove the following.
\begin{theorem}\label{strictinclusion}
The inclusion $H^{1}\subset H^{1}_{\rm{Riesz}}$ is strict.
\end{theorem}
Thus in our setting there is no characterization of the atomic Hardy space by means
of Riesz transforms, in contrast to the situation in the Euclidean and many other cases,
 as we explain below.


In this paper, we also prove the following $H^1$-$H^1$ boundedness result.
\begin{theorem}\label{t: H1bdd}
The Riesz transforms $\mathcal R_i$, $i=1,2$, are bounded from $
H^{1} $ to $H^1$.
\end{theorem}
This theorem is analogous to other boundedness results for Riesz transforms in different contexts known in the literature, but its proof is very different since the Hardy space $H^1$ has only an atomic definition in our setting. Thus the proof requires an explicit construction of the atomic decomposition of $\mathcal R_iA$ for an atom $A$. This will be based on a delicate argument of ``mass transport" given in Section \ref{RiH1H1}.  It is still an open problem whether also $\mathcal R_0$ is bounded on $H^1$.

\bigskip

The relation between Hardy spaces and Riesz transforms has been studied in different settings in the literature. Here we shall only mention papers where either a characterization of an atomic Hardy space in terms of Riesz transforms or a boundedness result for Riesz transforms on an atomic Hardy space has been investigated.

%

In the Euclidean case, the atomic Hardy space $H^1(\mathbb R^n)$ can also be characterized by means of the first-order Riesz transforms, with equivalence of norms (see \cite[Chapter 3]{St}). 
Moreover, all Riesz transforms are bounded from $H^1(\mathbb R^n)$ to $H^1(\mathbb R^n)$. Analogous results have been proved on nilpotent Lie groups \cite{LoVa} and more generally on Lie groups of polynomial growth \cite{SC} for the first-order Riesz transforms associated with a sub-Laplacian and for the atomic Hardy space defined in this setting (see \cite{CW2}).
See \cite{AMR, CMM, MaRu, Ru} and the references therein for boundedness results on Hardy spaces for first-order Riesz transforms on various doubling Riemannian manifolds. The relation between Riesz transforms and Hardy spaces associated with particular classes of operators has also been studied,  for instance in \cite{BDG, DP, D, MMS}.

To the best of our knowledge, there are no results of this kind in the literature  for Hardy-type spaces on nondoubling Lie groups. Results on the $H^1$-$L^1$ boundedness of Riesz transforms on some Lie groups of exponential growth can be found in \cite{MaVa} and  \cite{SV1}, but it seems that the $H^1$-$H^1$ boundedness of such operators has not been investigated. Some results on the $H^1$-$L^1$ boundedness of Riesz transforms for a flow Laplacian on infinite trees equipped with flow measures of exponential growth are obtained in \cite{LMSTV, MSTV, MSV}; these can be thought as a discrete counterpart of the results on Lie groups mentioned above.

\smallskip

Finally we mention some negative results which recall in some sense the negative result in our Theorem \ref{strictinclusion}.

On homogeneous trees, \cite{CeMe} says that the Hardy space defined in terms of the first-order Riesz transform associated with the standard combinatorial Laplacian does not admit an atomic decomposition. Santagati recently showed in \cite{S} that the natural atomic Hardy space defined on homogeneous trees equipped with the canonical flow measure and the flow Laplacian does not admit a characterization in terms of the Riesz transform, proving a discrete counterpart of Theorem \ref{strictinclusion}.

In the setting of Riemannian manifolds of exponential growth with bounded geometry and spectral gap, it was proved in \cite{MMV} that the first-order Riesz transform is not bounded from the atomic Hardy space introduced in \cite{CaMaMe} into the space of integrable functions; hence it does not provide a characterization of this Hardy space. However, on the same class of manifolds, the Hardy space defined by means of the first-order Riesz transform was recently characterized as a suitably modified Hardy space of Goldberg type \cite{MaMeV, MaMeVaVe, MeVe}.

\bigskip

Our paper is organized as follows.  Section \ref{kernels} contains explicit formulas for the convolution kernels of the first-order Riesz transforms $\mathcal R_i$ and some of their derivatives. In Section \ref{strict}, we prove Theorem \ref{strictinclusion}, and Theorem \ref{t: H1bdd} is proved in Section~\ref{RiH1H1}.

\bigskip

In this paper, $C$ denotes a positive, finite constant which may vary from occurrence  to occurrence and may depend on parameters according to the context. Given two positive quantities $f$ and $g$, we mean by $f\lesssim g$ that there exists a constant $C$ such that $f\leq C\, g$, and $f\sim g$ means that $g\lesssim f\lesssim g$.

\section{The convolution kernels of the Riesz transforms}\label{kernels}
In this section, we write the formulas for the convolution kernels of the Riesz transforms of the first order and some of their derivatives, which were computed in \cite{SV1}. First recall that the convolution of two (suitable) functions $f,\,g$ on $G$ is
\begin{equation}\label{conv1}
f\ast g(x)=\int_Gf(xy^{-1})\,g(y)\dir (y)\qquad\forall x\in G\,.
\end{equation}
The convolution kernel from the right of the operator $\mathcal R_i=X_i\,\Delta^{-1/2}$ is the distribution ${\rm{pv~}} k_i$, where $k_i=X_iU$, and $U$ is the convolution kernel of $\Delta^{-1/2}$ given by
\begin{equation}\label{formulaU}
U(x)=\frac{1}{2\pi^{2}}\,\delta^{1/2}(x)\,\frac{1}{r(x)\,\sinh r(x)} \,.
\end{equation}
The integral kernel of $\mathcal R_i$ is given by the function
\begin{equation}\label{kernelRi}
R_i(x,y)=\delta(y)\,k_i(y^{-1}x)\qquad   x\neq y\,,
\end{equation}
in the sense that
\begin{equation*}
 \mathcal R_i f(x) =  f\ast  {\rm{pv~}}  k_i (x)  ={\rm{pv~}}  \int R_i(x,y) f(y)\,d\rho(y),  \qquad f \in C_0^\infty(G).
\end{equation*}

The explicit formulas for $k_i$ were obtained in \cite[formulae (2.7), (2.8)]{SV1} by means of \eqref{formulaU} and the fact that for $x \neq e$
\begin{equation}\label{derivativedistance}
X_ir(x)=\begin{cases}
\frac{a-a^{-1}-a^{-1}(x_1^2+x_2^2)}{2\,\sinh r(x)}=\frac{a-\cosh r}{\sinh r}&{\rm{if}}~i=0\\
\frac{x_i}{\sinh r(x)}&{\rm{if}}~i=1,2\,,
\end{cases}
\end{equation}
where $r=r(x)=r(x_1,x_2,a)$ (see \cite[Lemma 2.1]{SV1}).

For $i=1,2$ and $x\neq e$
\begin{align}\label{k1}
2\pi^2 k_i(x)  &= - \,a^{-1}\,x_i\,\frac{\sinh r+r\,\cosh r}{r^2\,\sinh^3 r}\,
\end{align}
 and
\begin{equation}\label{k0}
\begin{aligned}
2\pi^2 k_0(x)
= \,- \frac{\sinh r+r\,\cosh r}{r^2\,\sinh^3 r} + a^{-1}\,\frac{\cosh r\sinh r+r}{r^2\,\sinh^3 r}\,,
\end{aligned}
\end{equation}
as verified in  \cite[formula (2.8)]{SV1}.
For all $x \neq e$ the derivative $X_1k_1(x)$ can be computed using  \eqref{k1} and \eqref{derivativedistance}:
\begin{multline} \label{derk1}
2\,\pi^{2}\,X_1k_1(x)\\
= a^{-1}\,\frac{x_1^2}{\sinh r}\,\frac{2r^2\cosh^2 r+r^2+2\sinh^2r+3r\sinh r\cosh r}{r^3\sinh^4 r}
-
    \frac{\sinh r+r\cosh r}{r^2\sinh^3r}\,.
\end{multline}
In particular
\begin{equation}\label{stimaderk1}
\begin{aligned}
|X_1k_1(x)|&\lesssim a^{-1}\,\frac{x_1^2}{r\sinh^3 r} +\frac{1}{r\sinh^2 r}\lesssim \frac{1}{r\sinh^2 r}\qquad \forall x\in B_1^c\,,
\end{aligned}
\end{equation}
 the last step since  $a^{-1}\,|x|^2 < 2 \cosh r(x) \simeq \sinh r$ for these $x$.


Instead of computing the derivative $X_1k_0(x)$ by means of formulas \eqref{k0} and \eqref{derivativedistance}, we estimate it for $x\in B_1^c$. We need only observe that when we differentiate a power of $\cosh r$ or $\sinh r$ with respect to $r$, the order of magnitude does not change for $r>1$, whereas a power of $r$ gets smaller when differentiated. From  \eqref{k0} we then get
\begin{equation*}
|X_1k_0(x)|\lesssim \frac{|x_1|}{r\,\sinh^3 r}+  a^{-1}\, \frac{|x_1|}{r^2\,\sinh^2 r}
+ a^{-1}\, \frac{|x_1|}{r\,\sinh^4 r}.
\end{equation*}
The last term here is no larger than the first term, since $a^{-1} < 2\cosh r \simeq \sinh r$,
and we conclude that
\begin{equation}\label{derk0stima}
|X_1k_0(x)|\lesssim \frac{|x_1|}{r\,\sinh^3 r}+  a^{-1}\, \frac{|x_1|}{r^2\,\sinh^2 r}\qquad \forall x\in B_1^c\,.
\end{equation}

In the sequel we shall repeatedly use the following integration formula (see for instance \cite[Lemma 1.3]{CGHM}): for any radial function $f$ such that $\delta^{1/2}f$ is integrable
\begin{equation}\label{intformula}
\int_G\delta^{1/2}f\dir=\int_0^{\infty}f(r)\,r\,\sinh r\di r\,.
\end{equation}

\section{Proof of Theorem \ref{strictinclusion} }\label{strict}

We recall a family of functions in the atomic Hardy space, which we introduced in \cite{SV2}. Let $L>2$
be large, and consider the rectangles $P_0=[-1,1]\times[-1,1]\times[\frac{1}{\nep},\nep]$ and
$P_L=(L,0,1)\cdot P_0=[L-1,L+1]\times[-1,1]\times[\frac{1}{\nep},\nep]$. We then define  $f_L=\chi_{P_L}-\chi_{P_0}$. Obviously $f_L$ is a multiple of an atom and it was proved in \cite{SV2} that $\|f_L\|_{H^1}\sim\,\log L$ for  $L>2$\,.
We claim that there exists a positive constant $C$ such that for $L>C$
\begin{equation}\label{Claim}
\|\mathcal  R_i f_L\|_{1}\lesssim \,\log\log L, \qquad i=0,1,2\,.
\end{equation}
Once this has been verified, we deduce that
$$
\lim_{L\rightarrow \infty}   \frac{ \|f_L\|_{H^1 }  }{\|f_L\|_{H^1_{\rm{Riesz}}}}=\infty\,,
$$
which proves Theorem \ref{strictinclusion}.

\smallskip

When we now prove \eqref{Claim}, we will neglect the case   $i=2$, since it is completely analogous to $i=1$.

\smallskip

Denote by $2P_0$ the rectangle $[-2,2]\times[-2,2]\times[\frac{1}{\nep^2},\nep^2]$ and by $2P_L $ the rectangle $(L,0,1)\cdot (2P_0)$.
For $i = 0$ and $1$
we shall estimate the $L^1$-norm of $\mathcal R_i f_L$  by integrating  over different regions, in four steps.

\bigskip

\noindent{\bf{Step 1.}} Since $\mathcal R_i$ is bounded  on $L^2$ and $\rho(2P_0)=\rho(2P_L)   \sim 1$, we get by applying the Cauchy--Schwarz inequality
\begin{equation}\label{int2R02R1}
\|\mathcal R_i f_L\|_{L^1(2P_0\cup 2P_L)}\lesssim \rho(2P_0\cup 2P_L)^{1/2}\,\|\mathcal R_i\|_{2\rightarrow 2}\,\|f_L\|_2
\simeq 1\,.
\end{equation}

\bigskip

\noindent{\bf{Step 2.}} Choose a ball $B=B(e,r_B)$ with $r_B=(\log L)^{2}$.
Then \eqref{metrica} implies that $B\supset 2P_0\cup 2P_L$ if $L$ is large enough. For any $x$
\begin{equation}
|\mathcal R_i \chi_{P_0}(x)|= \Big|\int \chi_{P_0}(y)\,k_i(y^{-1}x)\dil (y)\Big|\lesssim
\sup_{y \in P_0}  |k_i(y^{-1}x)| \,,
\end{equation}
since  $\lambda(P_0) \sim 1$.
If $x=(x_1,x_2,a)\in (2P_0)^c$ and $y=(y_1,y_2,b)\in P_0$, then
$$
|r(y^{-1}x)-r(x)|=|d(y,x)-d(x,e)|\leq d(y,e)\leq C\,.
$$
We will often use the following simple quotient formula:
\begin{equation}\label{y-1x}
y^{-1}x=(y_1,y_2,b)^{-1}(x_1,x_2,a)   =  ( b^{-1}(x_1-y_1), b^{-1}(x_2-y_2),ab^{-1} )    \,.
\end{equation}
Since here $ab^{-1}\sim a$ and $b^{-1}|x_1-y_1|\lesssim \max(1,|x_1|)$\,, applying
\eqref{k1} and \eqref{k0} we obtain for   $x\in (2P_0)^c$
\begin{equation}\label{R1chiR0}
|\mathcal R_1 \chi_{P_0}(x)|\lesssim \sup_{y \in P_0} |k_1(y^{-1}x)|\lesssim  a^{-1}\,  \frac{1+|x|}{r(x)\,\sinh^2 r(x)} \,
\end{equation}
and
\begin{equation}\label{R0chiR0'}
|\mathcal R_0 \chi_{P_0}(x)|\lesssim \sup_{y \in P_0} |k_0(y^{-1}x)|\lesssim
\frac{1}{r(x)\sinh^2r(x)}+ a^{-1}\,\frac{1}{r(x)^2\sinh r(x)}\,.
\end{equation}
We shall now integrate these quantities over the set $B\setminus {2P_0}$. For a point $x\notin 2P_0$ it is elementary to verify that \eqref{metrica} implies $r(x)>1$. If also $x=(x_1,x_2,a)\in B$, the same formula shows that $e^{-r_B} \le a \le  e^{r_B}$ and also
$$
a+a^{-1}\lesssim \cosh r(x)\simeq \sinh r(x).
$$
 For the expression in \eqref{R1chiR0} we get
\begin{align*}
\int_{B\setminus 2P_0} a^{-1}\,\frac{1+|x|}{r(x)\sinh^2 r(x)}\,\dir(x)
\lesssim & \int_{e^{-r_B}}^{e^{r_B}} \frac{da}{\log(a+a^{-1})}\,
\int \frac{1+|x|}{(a^2 +  1+|x|^2)^2}\,dx \\
\simeq &\, \int_{e^{-r_B}}^{e^{r_B}} \frac{da}{(a^2+1)\log(a+a^{-1})}\\
\lesssim & 1.
\end{align*}
The argument for the first summand in \eqref{R0chiR0'} is similar and left to the reader. For the second summand in  \eqref{R0chiR0'}, we use the integration formula \eqref{intformula} and the fact that
$B(e,1) \subset 2P_0$.  This gives
\begin{align*}
\int_{B\setminus 2P_0} a^{-1}\,\frac{1}{r(x)^2\sinh r(x)}\,\dir(x)
\lesssim & \int_1^{{r_B}} \frac{r\sinh r}{r^2 \sinh r}\,dr
\lesssim  \log r_B \simeq  \log \log L.
\end{align*}

 We conclude that for $i=1$ and $i=0$
 \begin{equation}       \label{intBBBB}
\int_{B\setminus 2P_0} |\mathcal R_i \chi_{P_0}|\dir\lesssim \log\log L\,.
\end{equation}
Here we observe that this inequality remains true, with the same proof, if $B$ is replaced by
the doubled ball $2B = B(e,2r_B)$.

The inequality \eqref{intBBBB} implies a similar estimate for $\mathcal R_i \chi_{P_L}$ in $B\setminus 2P_L$. Indeed,
 $\mathcal R_i\chi_{P_L}(x)=\mathcal R_i\chi_{P_0}(\tau_L x)$, where $\tau_Lx=(-L,0,1)\cdot x$, for any $x$ in $G$. Using the facts that $\tau_L P_L=P_0$ and $\tau_LB\subset 2B$ for large $L$, and changing variable $\tau_Lx=v$, we obtain
\begin{equation}\label{intBB}
\begin{aligned}
\int_{B\setminus 2P_L} |\mathcal R_i\chi_{P_L}(x)|\dir(x)&= \int_{B\setminus 2P_L} |\mathcal R_i\,\chi_{P_0}(\tau_Lx)|\dir(x)\\
&=\int_{\tau_LB\setminus \tau_L(2P_L)} |\mathcal R_i\chi_{P_0}(v)|\, \delta(-L,0,1) \dir(v)\\
&\leq \int_{2B\setminus 2P_0}|\mathcal R_i\chi_{P_0}|\,  \dir\\
&\lesssim \log\log L\,.
\end{aligned}
\end{equation}

From \eqref{intBBBB} and \eqref{intBB},  it now follows that
\begin{equation*}    
\int_{B\setminus (2P_0\cup\, 2P_L)} |\mathcal R_i f_L|\dir\lesssim \log\log L\,, \qquad i = 0, 1,
\end{equation*}
which ends Step 2.

\bigskip


\noindent{\bf{Step 3.}}
To deal with the complement of  the ball $B$, we will use  cancellation between the two parts of $f_L$.
  For a point $x\in B^c$, we write the convolution $f_L\ast k_i(x)$   as
\begin{align}\label{conv}
f_L\ast k_i(x)&=\
\int_{P_L}\,k_i\big(y^{-1} x\big)\di\lambda(y)-\int_{P_0}\,k_i(y^{-1}x)\di\lambda(y) \\
&=\int_{P_0}\,\big[k_i\big(y^{-1}(-L,0,1)x\big)-k_i(y^{-1}x)\big]\di\lambda(y)\,. \notag
\end{align}
Let now $y^{-1}=(y_1,y_2,b)$ be any point in $(P_0)^{-1}$ and $x=(x_1,x_2,a)$ any point in $B^c$. Then $y^{-1}(-L,0,1)x=y^{-1}(-L,0,1)y\cdot y^{-1}x=(-bL,0,1)\,y^{-1}x$, and the Mean Value Theorem implies
\begin{align} \label{meanvth}
k_i\big(y^{-1}(-L,0,1)x\big)-k_i(y^{-1}x) &=-bL\,\partial_1k_i\big((s,0,1)\,y^{-1}x\big),
\end{align}
for some $s\in (-bL,0)$.

 By the triangle inequality,
$$
\begin{aligned}
|r\big((s,0,1)\cdot y^{-1}x\big)-r(x)|&=|d\big(x,y\,\cdot(-s,0,1)\big)-d(x,e)|\\
&\leq d\big(y\,\cdot(-s,0,1),e\big)\\
&\leq d\big(y\,\cdot(-s,0,1),y \big)+d\big(y,e\big)\\
&=d\big((-s,0,1),e\big)+d\big(y,e\big)\\
&\lesssim \log L\,,
\end{aligned}
$$
where we also used \eqref{metrica}.  Now $r(x)>(\log L)^{2}$ implies $r\big((s,0,1)\,y^{-1}\,x\big)\sim r(x)$, and also $r\big((s,0,1)\,y^{-1}\,x\big) > 1$, since $L$ is large. Further,
\begin{equation}\label{sh}
\sinh r\big((s,0,1)\cdot y^{-1}x\big)\gtrsim  \sinh \left(r(x)-\log L\right)\gtrsim  \frac{\sinh r(x)}{L}\,.
\end{equation}

In the rest of  this step, we will deal only with the case $i=1$.

The estimate  \eqref{stimaderk1} now  implies, since $X_1=a\partial_1$ and $b \simeq 1$,
\begin{equation*}
\big|\partial_1k_1\big( (s,0,1)y^{-1}x    \big)  \big|\lesssim
(ba)^{-1}\,
\frac{1}{r\big((s,0,1)\,y^{-1}\,x \big)\sinh ^2 r\big((s,0,1)\,y^{-1}\,x \big)}\,.
\end{equation*}
 We deduce from \eqref{meanvth} and \eqref{sh} that for $x \in B^c$ and $y \in P_0$
\begin{align}\label{suptdifferenza}
\big| k_1\big(y^{-1}(-L,0,1)x\big)-k_1(y^{-1}x) \big| & \lesssim L\,a^{-1}\,\frac{1}{r\big((s,0,1)\,y^{-1}\,x \big)\,\sinh ^2 r\big((s,0,1)\,y^{-1}\,x \big)} \\ & \lesssim L^3 \,a^{-1}\,\frac{1}{r(x)  \sinh ^2 r(x)}. \notag
\end{align}

The same bound holds for $|\mathcal R_1 f_L|$  in $B^c$ because of  \eqref{conv}.

By applying the integration formula \eqref{intformula}, we finish Step 3 concluding that
\begin{equation}\label{intOmega2}
\begin{aligned}
\int_{B^c}|\mathcal R_1f_L |\dir \lesssim L^3\int_{r_B}^{\infty}\frac{1}{r\,\sinh^2 r}\,r\,\sinh r\,\di r
\lesssim  \frac{L^3}{\nep^{r_B}}\lesssim 1\,.
\end{aligned}
\end{equation}

\bigskip

\noindent{\bf{Step 4.}} It remains to take
 $i=0$ and estimate the integral of $\mathcal R_0 f_L$ over the complement of $B.$ This
 requires some modifications from the argument in the preceding step.

We split $B^c$ into  the following three parts:
$$
\begin{aligned}
\Gamma_1&=\left\{x=(x_1,x_2,a)\in B^c:\,a<a^*\quad \text{and} \quad |x|<f(a)\right\}\,,\\
\Gamma_2&=\left\{x=(x_1,x_2,a)\in B^c:\,a\geq a^* \right\}\,,
    \\ \Gamma_3&=\Big\{x=(x_1,x_2,a)\in B^c:\,a<a^*,\;\;|x|\geq f(a)\Big\}\,,
\end{aligned}
$$
where $a^*=\nep^{   -r_B/8}$                       and $f(a)=\exp \sqrt{\log a^{-1}}$\,.

To  integrate  $\mathcal R_0 f_L$ over  $\Gamma_1$, we consider
$\mathcal R_0\chi_{P_0}$ and $\mathcal R_0\chi_{P_L}$ separately, and apply \eqref{R0chiR0'}.
Observe that any point $x \in \Gamma_1$ satisfies
$\sinh r(x) \simeq \cosh r(x) > a^{-1}(1+|x|^2)/2$ and thus $r(x) \gtrsim \log a^{-1} $.
 Starting with $\mathcal R_0\chi_{P_0}$ and the last term in  \eqref{R0chiR0'}, we have
\begin{align*}
\int_{\Gamma_1} a^{-1}\,\frac{1}{r(x)^2\sinh r(x)}\,\dir(x)
 &\lesssim \int_0^{a^*}\frac{\di a}{a}\, a^{-1}\int_{|x|<f(a)} \frac{\di x}
 {\left(\log a^{-1}\right)^2\,a^{-1}\,(1+|x|^2)}  \\
 & \simeq \int_0^{a^*}\frac{\di a}{a} \frac{\log f(a)}  {\left(\log a^{-1}\right)^{2}  } \\
 & = \int_0^{a^*}\frac{\di a}{a \left(\log a^{-1}\right)^{3/2}  } \lesssim 1.
\end{align*}
  The first term in   \eqref{R0chiR0'} can be treated in a similar but easier way, and so
\begin{align*}
  \int_{\Gamma_1} \mathcal R_0\chi_{P_0}\,\dir \lesssim 1.
\end{align*}
  A translation argument like \eqref{intBB} shows that the same estimate holds for $\mathcal R_0\chi_{P_L}$,
  and thus for $\mathcal R_0 f_L$.

In order to estimate the integrals over $\Gamma_2$ and $\Gamma_3$, we
will use \eqref{conv} and \eqref{meanvth}, with $i=0$.  Combining  \eqref{meanvth} and
 \eqref{derk0stima}, we see that
$$
\big|\partial_1k_0\big( (s,0,1)y^{-1}x\big)  \big|\lesssim
 (ba)^{-1}\,\frac{|s+y_1+bx_1|}{r\,\sinh^3 r}+  (ba)^{-2}\, \frac{|s+y_1+bx_1|}{r^2\,\sinh^2 r}\,,
$$
where $ r=r\big(      (s,0,1)\,y^{-1}\,x      \big)\gtrsim 1$\,.
Using \eqref{sh} and the facts that $b \simeq 1$ and $0< s \lesssim L$, we conclude that
\begin{equation}\label{diffk0}
\begin{aligned}
\big|    k_0\big(y^{-1}(-L,0,1)x\big)-k_0(y^{-1}x) \big|
\lesssim L^C \frac{1+|x|}{a\,r(x)\,\sinh^3r(x)}+ L^C\frac{1+|x| }{a^2\,r(x)^2\,\sinh^2r(x)}
\,,
\end{aligned}
\end{equation}
for some $C$.

To integrate the first summand here, we neglect the factor $r(x)$ and apply the following estimate,
valid  for $r(x)>r_B$,
\begin{align}\label{sinh}
\sinh r(x)  \simeq \cosh r(x) \gtrsim e^{r_B} + a +  a^{-1}(1+|x|^2).
\end{align}
 Further,  we extend the integration to $B^c$ which contains $\Gamma_2 \cup \Gamma_3$, and get
\begin{align*}
\int_{B^c} L^C\, \frac{1+|x|}{a\,r(x)\,\sinh^3r(x)}\, \dir(x)
\lesssim &\: L^C \,\int_0^\infty \frac{da}{a} \int \frac{1+|x|}{a\,\left(e^{r_B} + a +  a^{-1}(1+|x|^2)\right)^3}\,dx \\
=  & \:L^C \,\int_0^\infty a\,da  \int \frac{1+|x|}{\left(a\,e^{r_B} + a^2 +  1+|x|^2\right)^3}\,dx \\
\lesssim &\: L^C \, \int_0^\infty {a}\,\left(a\,e^{r_B} + a^2\right)^{-3/2} \,da \\
\lesssim &\: L^C \, \int_0^{e^{r_B}} a^{-1/2}\, e^{-3r_B/2} \,da + L^C \int_{e^{r_B}}^\infty a^{-2} \,da \\
\lesssim &\: L^C \,e^{-r_B}
\lesssim 1\,.
\end{align*}

To integrate the second term in \eqref{diffk0}, at first over  $\Gamma_2$, we use \eqref{sinh} and argue almost as above. Thus
\begin{align*}
\int_{\Gamma_2} L^C\, \frac{1+|x|}{a^2\,r(x)^2\,\sinh^2r(x)}\, \dir(x)
=  & \:L^C \,\int_{a^*}^\infty \frac{da}{a}   \int \frac{1+|x|}{\left(a\,e^{r_B} + a^2 +  1+|x|^2\right)^2}\,dx \\
\lesssim &\: L^C \, \int_{a^*}^\infty  \frac{1}{a\,\left(a\,e^{r_B} + a^2\right)^{1/2}}\,da\\
\lesssim &\: L^C \, \int_{a^*}^{e^{r_B}} a^{-3/2}\, e^{-r_B/2} \,da + L^C \, \int_{e^{r_B}}^\infty a^{-2} \,da \\
\lesssim &\: L^C \, {(a^*)}^{-1/2}\,e^{-r_B/2} + L^C \, e^{-r_B}
\lesssim 1.
\end{align*}

Integrating the same term over $\Gamma_3$ in a similar way, we get at most
\begin{align*}
   L^C \,\int_0^{a^*} \frac{da}{a}   \int_{|x|>f(a)} \frac{1+|x|}{\left(1+|x|^2\right)^2}\,dx
   \lesssim \: L^C \,\int_0^{a^*}  \frac{da}{a\,f(a)}
 = &\: L^C \,\int_0^{a^*}  \frac{da}{a\,\exp \sqrt{\log a^{-1}} } \\
 = &\: L^C \int_{\sqrt{\log a^{-1}}}^\infty \, \frac{2s}{e^s}\,ds
   \\ \simeq &\: L^C \,  \frac{\sqrt{\log \,(a^*)^{-1}}}{\exp \sqrt{\log \,(a^*)^{-1}}}
\lesssim 1\,;
\end{align*}
here the change of variables used is $s = \sqrt{\log a^{-1}}$\,.
\vskip4pt

Summing up Step 4, we conclude
\begin{equation*}
\begin{aligned}
\int_{B^c}|\mathcal R_0f_L |\dir \lesssim 1\,.
\end{aligned}
\end{equation*}

\bigskip

The combined results of Steps 1--4 now prove
 \eqref{Claim} and thus also Theorem \ref{strictinclusion}.


\section{Proof of Theorem \ref{t: H1bdd}}\label{RiH1H1}

We give the proof only for  $i=1$, since the case  $i=2$ is completely
analogous.
It is enough to show that
\begin{equation}\label{uniformboundednessatom}
\sup\{\|\mathcal R_1A\|_{H^1}:A\,\,\mathrm{is \,\,a}\;  (1,\infty)\mathrm{-atom}\}<\infty\,.
\end{equation}

Indeed, let $f\in H^1$ and write $f=\sum_j \lambda_j A_j$, where $A_j$ are $(1,\infty)$-atoms and
$\sum_j| \lambda_j | <2\|f\|_{H^1}$. Since by \cite[Theorem 2.3]{HS} $\mathcal R_1$ is of weak type $(1,1)$ and the sum $\sum_j\lambda_jA_j$  is convergent in $L^1$,
$$
\mathcal R_1 f=\sum_j\lambda_j \mathcal R_1A_j
$$
with convergence in $L^{1,\infty}$, so that \eqref{uniformboundednessatom} would imply
$$
\|\mathcal R_1f\|_{H^1}\leq C \sum_j | \lambda_j |\leq C\,\|f\|_{H^1}\,
$$
and thus Theorem  \ref{t: H1bdd}.

To prove \eqref{uniformboundednessatom} we take a $(1,\infty)$-atom $A$ supported in a  \CZ set $R= [ -L/2,L/2]^2\times [\nep^{-\alpha},\nep^{\alpha}]$ centred at $(0,0,1)$. We call $A$ a large atom if $\alpha\geq 20$, a small atom if $\alpha<1$ and an intermediate atom if $1\leq \alpha<20$.

For each $x\in G$, we have
\begin{equation}\label{f1}
|\mathcal R_1A(x)|\leq \int_R|A(y)|\,|k_1(y^{-1}x)|\delta(y)\di \rho(y)\leq \rho( R )^{-1} \int_R |k_1(y^{-1}x)|\delta(y)\di \rho(y)\,.
\end{equation}
With  $y \in R$, we will use \eqref{k1} to estimate the value of $k_1$ at the point $
z=y^{-1}x$, in a way depending on the position of $x$ in the group $G$. We first observe that  \eqref{metrica} implies for a point $z=(z_1,z_2,d)$ such that $d<1/2$
$$
\sinh  r(z)\sim \cosh r(z)\sim d^{-1}\,(1+|z|^2)
$$
and $r(z)\gtrsim 1$.
From \eqref{k1} we then get
\begin{equation}\label{d-piccolo}
|k_1(z)|\lesssim d^{-1}\,|z|\,\frac{1}{r(z)\,\sinh^2 r(z)}\lesssim \frac{d\,|z|}{(1+|z|^2)^2}\,.
\end{equation}
On the other hand, if $d>2$ then
\begin{equation}\label{dlarge}
  \cosh r(z) =\frac{d+d^{-1}|z|^2}{2}\, \Big(   1+ O\Big( \frac{1}{d^2+|z|^2} \Big)  \Big)\,,
\end{equation}
and the same equality holds for  $\sinh r(z)$.
Further,
\begin{equation}
  \label{logr}
  r(z)\sim \log \frac{d^{2}+|z|^2}d \gtrsim \log d.
\end{equation}
 From \eqref{k1} we conclude when $d>2$
\begin{equation}\label{cgrande}
\begin{aligned}
2\pi^2 k_1(z)&  =-d^{-1}z_1 \frac{ \frac{d+d^{-1}|z|^2}{2}   \Big(   1+ O\Big( \frac{1}{d^2+|z|^2} \Big)  \Big)+r(z)\frac{d+d^{-1}|z|^2}{2}   \Big(   1+ O\Big( \frac{1}{d^2+|z|^2} \Big)  }{
r^2(z) \Big(\frac{d+d^{-1}|z|^2}{2}\Big)^3   \Big(   1+ O\Big( \frac{1}{d^2+|z|^2} \Big) \Big) }\\
&= - \frac{     4 dz_1      }{    (d^2+|z|^2)^2    }\frac{1+r(z)}{r^2(z)}\Big(   1+ O\Big( \frac{1}{d^2+|z|^2} \Big)\Big)\,.
\end{aligned}
\end{equation}

Since $A$ has vanishing integral, so has $\mathcal R_1A$.
Moreover, we shall see that  $\mathcal R_1A$ is integrable and has vanishing integral on each horizontal plane  $\Bbb R^2 \times \{a \}$ for  $a< \nep^{-\alpha}$ and for
 $a >\nep^{\alpha}$. We first verify that the function
$x \mapsto k_1(y^{-1}x)$ has these properties for any $y= (y_1, y_2, b) \in R$.

If $x = (x_1, x_2, a)$ with
 $a \notin [\nep^{-\alpha},\nep^{\alpha}]$ and $y \in R$, we see from \eqref{y-1x}
that $r(y^{-1}x) \gtrsim 1.$ Then \eqref{k1} implies that
\[
\int |k_1(y^{-1}(x_1, x_2,a)|\,\di x_1\di x_2 < \infty\,.
\]
 Now  $r(y^{-1}x)$ is a function of $a$, $b$ and $(x_1-y_1)^2 + (x_2-y_2)^2$,
and from  \eqref{k1} we see that $k_1(y^{-1}(x_1, x_2,a))$ is an odd function
in  $x_1-y_1$. Thus
\begin{equation}
  \label{simmetria}
\int k_1(y^{-1}(x_1, x_2,a))\,\di x_1\di x_2 = 0.
\end{equation}
Integrating this equality against $A(y)\delta(y)\,\di\rho(y)$, we conclude that
\begin{equation}\label{vanint}
  \int\mathcal R_1A(x_1, x_2,a) \,\di x_1\di x_2 = 0,  \qquad
a \notin [\nep^{-\alpha},\nep^{\alpha}],
\end{equation}
and the integral is absolutely convergent.

\smallskip

We treat the three atom sizes separately.

\bigskip

\noindent{\bf{Case I: large atom\,.}} In this case $\nep^{2\alpha}\leq L<\nep^{8\alpha}$. We will construct an atomic decomposition of $\mathcal R_1 A$.

\smallskip

Let us first give a rough description of the idea. We will split
 $G$ into slices of type $S_k = \Bbb R^2 \times I_k$, for disjoint intervals $I_k$, $k\in\mathbb Z$. Then \eqref{vanint} will imply that
\begin{equation}\label{vanint2}
  \int_{S_k} \mathcal R_1A \,\di\rho = 0,
\end{equation}
which will make it possible to decompose  $\mathcal R_1A\, \chi_{S_{k}}$ into atoms.

For each $k$, the slice $S_k$ will be split into disjoint sets $B_{k,j}=  \tilde B_{k,j} \times I_k$,~\;$j=0,1,\dots$, where the $\tilde B_{k,j}$,~$j=0,1,\dots$, form a partition of $\mathbb R^2$ and expand exponentially towards infinity as $j\rightarrow +\infty$. We want to make each $\mathcal R_1A\, \chi_{B_{k,j}}$ into an atom multiple,
supported
in a suitable  Calder\'on-Zygmund set  $Z_{k,j} \supseteq B_{k,j} $.
The $Z_{k,j} $
will not be disjoint but increasing in $j$. Then
\[
\mathcal R_1A\, \chi_{B_{k,j}} -  \int_{B_{k,j}} \mathcal R_1A \,\di \rho \,\,
 \frac{\chi_{Z_{k,j}}}{\rho(Z_{k,j})}
\]
 has moment 0 and is an atom multiple. But these functions do not sum up in $j$
to $\mathcal R_1A\, \chi_{S_{k}}$ as we want. Instead, we will modify
$\mathcal R_1A\, \chi_{B_{k,j}}$ by two quantities involving
$B_{k,\ell}$ for all $\ell \ge j$. This can be seen as a way of transporting the ``mass" $\int_{B_{k,0}}\mathcal R_1A \,\di \rho $ in the innermost $B_{k,0}$
step by step through $B_{k,1}, B_{k,2}, \dots$ towards infinity. At each step one leaves in $B_{k,j}$ the mass $- \int_{B_{k,j}} \mathcal R_1A \,\di \rho$ needed there. This leads to a telescopic sum,
and produces atom multiples summing up to $\mathcal R_1A\, \chi_{S_{k}}$.

\smallskip

Let us now go into the details of this construction. We decompose $G$ as the union of the following three regions:
\begin{align*}
\Omega_1&=\{x=(x_1,x_2,a)\in G: \,a <\nep^{-1-3\alpha/2} \}\,,\\
\Omega_2&=
\{x=(x_1,x_2,a)\in G: \,\nep^{-1-3\alpha/2}<a<\nep^{1+3\alpha/2} \}\,,\\
\Omega_3&=\{x=(x_1,x_2,a)\in G: \,a>\nep^{1+3\alpha/2} \}\,.
\end{align*}
For a point $x$ in any of these regions,
we will derive size estimates  of $\mathcal R_1A(x)$ via estimates of the
kernel $k_1(  y^{-1}x)$, where  $y=(y_1,y_2,b)\in R$. Observe  that
if $|x|>2L$ then $L < |x-y|\sim |x|$, but if  $|x|\le 2L$ then
 $|x-y|\le 3L$.
\smallskip

Let $x\in \Omega_1$, so that $b^{-1}a<1/4$.
Consider first the case $|x|>2L$. With $z=y^{-1}x$ in \eqref{d-piccolo}, we get
\[
|k_1(  y^{-1}x)|\lesssim   \frac{b^{-1}a\,b^{-1}|x-y|}{(1+b^{-2}|x-y|^2)^2}
\lesssim \frac{a\,b^{2}}{|x|^3}\,.
\]
We conclude from \eqref{f1} that for  $|x|>2L$
  \begin{equation}  \label{Omega1}
\begin{aligned} 
|\mathcal R_1A(x)|\lesssim \frac{1}{\alpha L^2}  \int_{\nep^{-\alpha}}^{\nep^{\alpha}}\int_{|y|<L/2}  \frac{ab^{2}}{|x|^3}        \,b^{-2}\, \frac{\di y\,\di b}{b}
\lesssim \frac{a}{|x|^3}\,.
\end{aligned}
  \end{equation}
In the case when  $|x|\le 2L$,  we similarly conclude from \eqref{d-piccolo} and \eqref{f1}
\begin{equation}\label{Omega2}
\begin{aligned}
|\mathcal R_1A(x)|&\lesssim  \frac{1}{\alpha L^2} \int_{\nep^{-\alpha}}^{\nep^{\alpha}}\int_{|y|<L/2}    \frac{ab^{-2}\,|x-y|}{(1+b^{-2}|x-y|^2)^2}  \,b^{-2} \frac{\di y\,\di b}{b}\\
&= \frac{a}{\alpha L^2}  \int_{\nep^{-\alpha}}^{\nep^{\alpha}}\int_{|v|<3L/b}     \frac{b^{-1}\,|v|}{(1+ |v|^2)^2}     \,  \frac{\di v  \,  \di b}{b}\\
&\lesssim \frac{a}{\alpha L^2} \int_{\nep^{-\alpha}}^{\nep^{\alpha}}
\frac{\di b}{b^2}\\
&\leq \frac{a \nep^{\alpha}}{\alpha L^2}\,.
 \end{aligned}
\end{equation}
Let for $k = -1,-2,\dots$
\begin{equation*}
  S_k = \Bbb R^2 \times  \left[\nep^{-2^{|k|}-3\alpha/2} , \nep^{-2^{|k|-1}-3\alpha/2}\right]
\end{equation*}
and
\begin{equation*}
  B_{k,j}= \{x\in S_k:\,  2\nep^{j-1} L \leq |x|\leq 2\nep^jL   \}
\end{equation*}
for $j = 1,2,\dots$, but for $j=0$ instead
\begin{equation*}
  B_{k,0}= \{x\in S_k:\, |x|\leq 2L   \}.
\end{equation*}
Except for boundaries, $S_k$ is the disjoint union of the $B_{k,j}$,
and  $\Omega_1$ is  the disjoint union of the $S_k$.
The measure of $B_{k,j}$ is  $\rho(B_{k,j}) \sim e^{2j} L^2 2^{|k|}$.
For $x \in  B_{k,j}$ with $j\ge 1$, \eqref{Omega1} implies that
\begin{equation}
  \label{kj}
  |\mathcal R_1A(x)| \lesssim \frac{\nep^{-2^{|k|-1}-3\alpha/2}}{e^{3j}L^3}
 \lesssim \frac{\nep^{-2^{|k|-1}-\alpha/2}}{e^{3j}L^2}.
\end{equation}
The last bound here holds also for $j=0$, since \eqref{Omega2} yields
\begin{equation*}
   |\mathcal R_1A(x)|
 \lesssim \frac{\nep^{-2^{|k|-1}-\alpha/2}}{L^2},
\end{equation*}
if  $x \in B_{k,0}$. For $k\le -1$ and $j\ge 0$, we define
\begin{equation*}
  Z_{k,j}=[-2 \nep^j L,2\nep^jL]^2\times [\nep^{-3\cdot 2^{|k|-1}-3\alpha-j/2} , \nep^{-2^{|k|-1}+j/2}]\,.
\end{equation*}
This set contains  $  B_{k,j}$ and is a Calder\'on--Zygmund set centered at
 $\left(0,0, \nep^{-2\cdot 2^{|k|-1}-3\alpha/2} \right)$, of parameter
 $ 2^{|k|-1}+3\alpha/2 +j/2$ and of measure
 $\rho(Z_{k,j})\sim \nep^{2j} \,L^2 \,(2^{|k|}+\alpha+j)$.

We define for $k\leq -1$ and $j = 0, 1,\dots$
 \begin{equation*}
   f_{k,j} = \rho(Z_{k,j})^{-1} \chi_{Z_{k,j}}
 \end{equation*}
and
\begin{equation*}
  m_{k,j} = \sum_{\ell = j}^\infty \int_{B_{k,\ell}} \mathcal R_1 A \,\di \rho.
\end{equation*}
Observe that \eqref{vanint} implies \eqref{vanint2} for each $k \le -1$, which
means that  $m_{k,0} =0$. For $j\ge 1$ we get from \eqref{kj}
\begin{equation}
  \label{mkj}
  |m_{k,j}| \lesssim \sum_{\ell = j}^\infty
\frac{\nep^{-2^{|k|-1}-\alpha/2}\,2^{|k|}}{\nep^{\ell}}
\sim \frac{\nep^{-2^{|k|-1}-\alpha/2}\,2^{|k|}}{\nep^{j}}.
\end{equation}
The functions
\begin{equation*}
  A_{k,j} = \mathcal R_1 A\, \chi_{B_{k,j}}
- m_{k,j}\,f_{k,j} +m_{k,j+1}\,f_{k,j+1}
\end{equation*}
 are bounded and supported in   $  Z_{k,j+1}$,
and they are easily seen to have  integrals 0. Thus they are multiples of
$(1,\infty)$-atoms.
Their sum over $j$ is
\begin{equation*}
  \sum_{j=0}^\infty  A_{k,j} = \mathcal R_1A\, \chi_{S_{k}}- m_{k,0}\,f_{k,0}
=  \mathcal R_1A\, \chi_{S_{k}}.
\end{equation*}
 Thus we have an atomic decomposition of
$\mathcal R_1A\, \chi_{S_{k}}$.
 Summing over  $k = -1,-2,\dots$, we get a decomposition
of $\mathcal R_1A\, \chi_{\Omega_1}$.
Combining  \eqref{kj} and  \eqref{mkj}, and multiplying   \eqref{kj}
 by  $\rho(Z_{k,j+1})$, we obtain
\begin{equation*}
  \| A_{k,j} \|_{H^1} \lesssim
\frac{\nep^{-2^{|k|-1}}\,\nep^{-\alpha/2}\, (2^{|k|}+\alpha+j+1)}{\nep^{j}}
+  \frac{\nep^{-2^{|k|-1}}\nep^{-\alpha/2}\,2^{|k|}}{\nep^{j}}.
\end{equation*}
The right-hand side here is summable in $j$ and $k$, and we conclude that
\begin{equation}
  \label{omega12}
 \| \mathcal R_1 A\, \chi_{\Omega_1}\|_{H^1} \lesssim 1.
\end{equation}

\bigskip

We next consider the region $\Omega_3 $, and estimate
$\mathcal R_1A$ there. Let $x\in\Omega_3$ with $|x|>2L$, and take $y=(y_1,y_2,b)\in R$. Setting $\tilde y=(0,0,b)$ we
decompose $\mathcal R_1A$   as
$$
\begin{aligned}
\mathcal R_1A(x)&=\int_RA(y)[k_1(y^{-1}x)\delta(y)-  k_1(\tilde y^{-1}x)\delta(y)]d\rho(y)\\
&+\int_RA(y)[k_1(\tilde y^{-1}x)\delta(y) -k_1(x) ]d\rho(y)\\
&=\mathcal R_1'A(x)+\mathcal R_1''A(x)\,.
\end{aligned}
$$
In this case $b^{-1}a>2$, and from \eqref{dlarge} we get for any $z=(z_1,z_2,b)$ with $|z_i|<L$
$$
\cosh r(y^{-1}x)\sim \cosh r(z^{-1}x)\sim b^{-1}a+ba^{-1}b^{-2}|x-z|^2\sim\frac{a}{b}+\frac{|x|^2}{ab}\,,
$$
and also
$$
r(y^{-1}x)\sim  r(z^{-1}x)\sim  \log\Big(\frac{a}{b}+\frac{|x|^2}{ab}\Big)
\ge \log a -   \log b \sim  \log a \,,
$$
since $\log a > 3\alpha/2$ and  $\log b <\alpha$.
By the Mean Value Theorem
$$
|k_1(y^{-1}x) -  k_1(\tilde y^{-1}x) |\leq L\,b^{-1} (|\partial_1k_1(z^{-1}x)| + |\partial_2k_1(z^{-1}x)|)\,,
$$
for some $z=(z_1,z_2,b)$ with $|z_i|<L$. Now apply the estimate \eqref{stimaderk1} for $X_1k_1=a\partial_1k_1$ to obtain that
\begin{equation}\label{stimader1k1}
 \left|   {\partial_1}  k_1(x)\right| \lesssim \frac1{ar(x)\cosh^2r(x)} \qquad {\rm{for\,\,}}\;\; r(x)>1\,.
\end{equation}
The same estimate holds for $|\partial_2k_1|$. We then have
$$
\left|k_1(y^{-1}x) -  k_1(\tilde y^{-1}x) \right|\,\lesssim\, L\,b^{-1} \,\frac{1}{b^{-1}a  \Big(\frac{a}{b}+\frac{|x|^2}{ab}\Big)^2 \, \log a}
= \frac{Lab^2}{(a^2+|x|^2)^2 \log{a}  }\,,
$$
so for $|x|>2L$ and $x\in \Omega_3$
\begin{equation}\label{R_1A'}
|\mathcal R_1'A(x)| \,\lesssim \, \frac{1}{\alpha L^2}\int_{\nep^{-\alpha}}^{\nep^{\alpha}}\int_{|y|<L/2}   \frac{Lab^2}{(a^2+|x|^2)^2 \log {a} }\,b^{-2}   \frac{\di y\di b}{b}
\,\leq\, \frac{La}{(a^2+|x|^2)^2 \log {a} }\,.
\end{equation}
For $\mathcal R_1'' A$ we use \eqref{cgrande} to see that
$$
2\pi^2\, k_1(\tilde y^{-1}x)\delta(y)=
-\frac{     4 ax_1   b^2   }{    (a^2+ |x|^2)^2    }\,
\frac{1+r( \tilde y^{-1}x  )}{r^2(\tilde y^{-1}x)}\,
\Big(   1+ O\Big( \frac{b^{2}}{a^2+|x|^2} \Big) \Big)\,b^{-2}.
$$
This must be compared with $2\pi^2\, k_1(x)$, which is obtained by replacing
$b$ by 1 and thus $\tilde y$ by $e$ here. Observe that

$$
\Big|\frac{1+r( \tilde y^{-1}x  )}{r^2(\tilde y^{-1}x)}-\frac{1+r( x  )}{r^2(x)}\Big|\lesssim \frac{|r( \tilde y^{-1}x  )-r(x)|}{r(\tilde y^{-1}x)r(x)} \lesssim \frac{  |\log b| }{(\log a)^2 }\,,
$$
where we applied the triangle inequality and \eqref{logr}.
Using this and estimating the $O(\dots)$ terms, we obtain
$$
\begin{aligned}
\left| k_1(\tilde y^{-1}x)\delta(y)-k_1(x) \right|
&\lesssim \frac{ a |x_1| }{(a^2+ |x|^2)^2}\frac{1+  |\log b| }{(\log a)^2 }\,.
\end{aligned}
$$
We then deduce that for $|x|>2L$
\begin{equation}\label{R_1A''}
  \begin{aligned}
|\mathcal R_1''A(x)|&\lesssim \frac{1}{\alpha L^2}\,\frac{     a |x_1|      }{    (a^2+ |x|^2)^2    }\frac{ 1 }{(\log a)^2 }\int_{|y|<L/2}\int_{\nep^{-\alpha}}^{\nep^{\alpha}} (1+|\log b|) \frac{\di b \di y}{b}\\
&\lesssim \frac{  \alpha   a |x_1|      }{    (a^2+ |x|^2)^2 \, (\log a)^2    }\,.
\end{aligned}
\end{equation}

Take now $x\in\Omega_3$ with $|x| \le 2L$. If also $y\in R$,  \eqref{cgrande} implies
$$
\begin{aligned}
|k_1(y^{-1}x)|\lesssim \frac{ a b^2|x-y|}{(a^2+|x-y|^2)^2   \,\log a}\,.
\end{aligned}
$$
For $\mathcal R_1A$ we then obtain
\begin{equation}\label{stimapuntualeOmega4}
\begin{aligned}
|\mathcal R_1A(x)|&\lesssim \frac{a}{\alpha L^2}\int_{\nep^{-\alpha}}^{\nep^{\alpha}}\int_{|y|<L/2}  \frac{    |x-y|}{(a^2+|x-y|^2)^2   \,\log a}\frac{\di y\di b}{b}\\
&\leq \frac{1}{\alpha L^2\log a}\int_{|v|<5L/2a}  \frac{    |v|}{(1+|v|^2)^2 }\,dv \,\int_{\nep^{-\alpha}}^{\nep^{\alpha}} \frac{\di b}{b}\\
&\lesssim \frac{1}{ L^2 \log a} \left(\min(1,L/a)\right)^3 \,.
\end{aligned}
\end{equation}
From this estimate we get for  $x\in \Omega_3$  and $|x| \le 2L$
\begin{equation}\label{agrandexpiccolo}
|\mathcal R_1A(x)|\lesssim \begin{cases}
\frac{L}{a^3\log a}&\mathrm{if}\,\,\; a>L\,,\\
\frac{1}{L^2\log a}&\mathrm{if}\,\,\; \nep^{1+\frac{3}{2}\alpha} \leq a\leq L\,.\\
\end{cases}
\end{equation}
 The atomic decomposition of $\mathcal R_1A$ in $\Omega_3$
is analogous to that in  $\Omega_1$.
But since the estimate \eqref{agrandexpiccolo} distinguishes between the cases $a>L$ and
 $a\le L$, we need to do the same when we define the slices  $S_k$.

For this, we first choose an integer $\kappa>1$ and a number $Q\in (2,3)$
such that $Q^\kappa = \log L - 3\alpha/2$. This is possible, since the equation
amounts to
\[
\kappa = \frac{\log (\log L - 3\alpha/2)}{\log Q}
\]
and $\log L - 3\alpha/2 \ge \alpha/2 \ge 10$. One can therefore choose  $\kappa$
between $({\log (\log L - 3\alpha/2)})/{\log 3}$ and $({\log (\log L - 3\alpha/2)})/{\log 2}$,  and then  $Q$ will be determined.

For $k= \kappa+1, \kappa+2, \dots$, we define slices
\begin{equation*}
  S_k = \Bbb R^2 \times [\nep^{k-\kappa-1}L, \nep^{k-\kappa}L].
\end{equation*}
But for  $k= 1,\dots \kappa$, the slices are instead
\begin{equation*}
  S_k = \Bbb R^2 \times [\nep^{1-Q^{\kappa-k+1}}L, \nep^{1-Q^{\kappa-k}}L].
\end{equation*}
It is easily verified that $(S_k)_1^\infty$ is a partition of
$\Omega_3$, except for boundaries, and that with $x\in S_k$ the cases $k > \kappa$
and $k \le \kappa$  correspond to  $a\ge L$ and $a\le L$, respectively.

For $k > \kappa$  we  define
\begin{equation*}
   B_{k,j}=
\{x\in S_k:\,  2\nep^{k-\kappa+j-1} L \leq |x|\leq 2\nep^{k-\kappa+j}L   \},
\qquad j = 1,2,\dots,
\end{equation*}
and
\begin{equation*}
   B_{k,0}= \{x\in S_k:\, |x|\leq 2\nep^{k-\kappa}L   \};
\end{equation*}
we also define
\begin{equation*}
   Z_{k,j}= [-2 \nep^{k-\kappa+j} L,2\nep^{k-\kappa+j}L]^2\times
[\nep^{k-\kappa-1-j/2}L , \nep^{k-\kappa+j/2}L]\,,
\qquad j = 0, 1,\dots.
\end{equation*}
The measures of these sets are given by
\begin{equation}\label{bigkmeas}
  \rho(B_{k,j}) \sim \nep^{2(k-\kappa+j)}L^2 \qquad \mathrm{and}
  \qquad \rho(Z_{k,j}) \sim  \nep^{2(k-\kappa+j)}L^2  (j+1),
\end{equation}
 and $Z_{k,j}$ is
a Calder\'on--Zygmund set, with center $(0,0, \nep^{k-\kappa-1/2}L)$ and parameter $(j+1)/2$.

For $k > \kappa$ and $j\ge 1$, we conclude from
\eqref{R_1A'} and \eqref{R_1A''} that for $x\in B_{k,j}$
\begin{equation}\label{bigk}
|\mathcal R_1A(x)|  \lesssim
\frac 1 {L^2\,\nep^{3(k-\kappa)}\,\nep^{4j}\,(k-\kappa+\alpha)} +
\frac \alpha  {L^2\,\nep^{2(k-\kappa)}\,\nep^{3j}\,(k-\kappa+\alpha)^2}.
\end{equation}
Observe that this holds also for  $j = 0$ because of
\eqref{agrandexpiccolo}.

\bigskip

For  $1 \le k \le \kappa$, we set instead
\begin{equation*}
   B_{k,j}=
\{x\in S_k:\,  2\nep^{j-1} L \leq |x|\leq 2\nep^{j}L   \},
\qquad j = 1,2,\dots,
\end{equation*}
and
\begin{equation*}
   B_{k,0}=
\{x\in S_k:\,  |x|\leq 2L   \}\,,
\end{equation*}
and also
\begin{equation*}
   Z_{k,j}= [-2 \nep^{j} L,2\nep^{j}L]^2\times
[\nep^{1-Q^{\kappa-k+1}-j/2}L , \nep^{1-Q^{\kappa-k}+j/2}L]\,,
\qquad j = 0, 1,\dots.
\end{equation*}
The measures are now given by
\begin{equation}\label{smallkmeas}
  \rho(B_{k,j}) \sim   \nep^{2j}L^2 Q^{\kappa-k} \qquad \mathrm{and}
  \qquad \rho(Z_{k,j}) \sim \nep^{2j}L^2 (Q^{\kappa-k}+j),
\end{equation}
 and $Z_{k,j}$ is again
a Calder\'on-Zygmund set, with center $(0,0, L\nep^{1-Q^{\kappa-k}(Q+1)/2})$ and parameter ${j}/{2}+Q^{\kappa-k}\,(Q-1)/2$.

If  $1 \le k \le \kappa$ and $j\ge 1$, we derive from \eqref{R_1A'},
\eqref{R_1A''}  and the inequality $\log a > \alpha$ that
\begin{equation}\label{smallk}
|\mathcal R_1A(x)|  \lesssim
\frac {\nep^{-Q^{\kappa-k}}} {L^2e^{3j}\,\alpha},
\end{equation}
for any $x\in B_{k,j}$. However, when   $x\in B_{k,0}$  we get from \eqref{agrandexpiccolo}
\begin{equation}
  \label{smallkj0}
  |\mathcal R_1A(x)|  \lesssim \frac 1
 {L^2\alpha}.
\end{equation}
For  $k\ge 1$ and  $j\ge 0$, we let as before
   $f_{k,j} = \rho(Z_{k,j})^{-1} \chi_{Z_{k,j}}$ and
\begin{equation*}
  m_{k,j} = \sum_{\ell = j}^\infty \int_{B_{k,\ell}} \mathcal R_1 A \,\di \rho.
\end{equation*}
Notice that $m_{k,0} = 0$ for all $k \ge 1$, because of \eqref{vanint}. From \eqref{bigk} and \eqref{bigkmeas}, we obtain for $k > \kappa$ and  $j\ge 0$
\begin{equation}\label{bigkm}
  |m_{k,j}| \lesssim \frac 1 {\nep^{k-\kappa}\,\nep^{2j}\,(k-\kappa+\alpha)} +
\frac \alpha  {\nep^{j}\,(k-\kappa+\alpha)^2}.
\end{equation}
For $1\le k \le \kappa$ and  $j\ge 1$, we similarly get
from \eqref{smallk} and  \eqref{smallkmeas}
\begin{equation}\label{smallkm}
  |m_{k,j}| \lesssim \frac {\nep^{-Q^{\kappa-k}}\,Q^{\kappa-k}} {\nep^{j}\,\alpha}.
\end{equation}
We now come to the atoms, and define for  $k \ge 1$ and  $j\ge 0$
\begin{equation} \label{atom}
  A_{k,j} = \mathcal R_1A\,\chi_{B_{k,j}}
-m_{k,j} f_{k,j} + m_{k,j+1} f_{k,j+1}.
\end{equation}
Then $ A_{k,j}$ is supported in the  Calder\'on--Zygmund set $Z_{k,j+1}$
and has vanishing moment. Further
\begin{equation}
  \label{total}
  \sum_{k=1}^\infty \sum_{j=0}^\infty A_{k,j}
= \mathcal R_1 A \,\chi_{\Omega_3}.
\end{equation}
The  $ A_{k,j}$ are multiples of  $(1,\infty)$-atoms.
We estimate their norms in $H^1$. When $k>\kappa$ and  $j\ge 0$,
\eqref{bigk}, \eqref{bigkmeas} and \eqref{bigkm} imply
\begin{align*}
 \| A_{k,j}\|_{H^1} \lesssim &\,
  \frac  {  j+1}  {\nep^{k-\kappa}\,\nep^{2j}\,(k-\kappa+\alpha)} +
\frac { (j+1)   \alpha}  {\nep^{j}\,(k-\kappa+\alpha)^2} \\
&+ \frac 1 {\nep^{k-\kappa}\,\nep^{2j}\,(k-\kappa+\alpha)} +
\frac \alpha  {\nep^{j}\,(k-\kappa+\alpha)^2}.
\end{align*}
Summing, one obtains
\begin{equation}
  \label{bigksum}
   \|\sum_{k>\kappa} \sum_{j\ge0} A_{k,j}\|_{H^1} \lesssim 1.
\end{equation}
For $1\le k \le \kappa$ and $j\ge 1$, we similarly see from
\eqref{smallk}, \eqref{smallkmeas} and \eqref{smallkm} that
\begin{align*}
 \| A_{k,j}\|_{H^1} \lesssim \,
\frac {\nep^{-Q^{\kappa-k}}\,(Q^{\kappa-k}+j) }  {\nep^{j}\alpha}
\end{align*}
Again, we can sum and get
\begin{equation}
  \label{smallksum}
   \|\sum_{1\le k\le \kappa} \sum_{j\ge1} A_{k,j}\|_{H^1} \lesssim 1.
\end{equation}
For $A_{k,0}$ with $ 1 \le k\le \kappa$, we have in view of
\eqref{smallkj0}, \eqref{smallkmeas} and \eqref{smallkm}
\begin{align*}
  \| A_{k,0}\|_{H^1} \lesssim  \frac {  Q^{\kappa-k}}
 {\alpha} +
 \frac {\nep^{-Q^{\kappa-k}}Q^{\kappa-k}} {\alpha},
\end{align*}
and since $Q^{\kappa} \simeq \alpha$, this implies
\begin{equation}
 \label{smallkj0sum}
   \|\sum_{1\le k\le \kappa} A_{k,0}\|_{H^1} \lesssim 1.
\end{equation}
Summing up, we have an atomic decomposition of
the restriction of $\mathcal R_1 A $ to $\Omega_3$ with control of the norm.

\smallskip

To deal with  the region $\Omega_2$, we first derive estimates for the
kernel and $\mathcal R_1A$ there. For $x = (x_1, x_2,a)\in\Omega_2$ with $|x| > 2L$ and $y\in R$, one has $|x-y|\sim |x|$ and
 $\nep^{-1-5\alpha/2}<b^{-1}a<\nep^{1+5\alpha/2}$, so that $|x|^2>4L^2>4\nep^{4\alpha}>2(a^2+b^2)$ and
$$
\cosh r(y^{-1}x)\sim \frac{a}{b}+\frac{b}{a}+\frac{b}{a}\,b^{-2}\,|x|^2\sim \frac{|x|^2}{ab}\qquad{\rm{and}}\qquad
r(y^{-1}x)\sim \log|x|\,.
$$
This implies that by \eqref{k1}
$$
|k_1(y^{-1}x)|\lesssim \frac{b}{a}\,|x|\,b^{-1}\Big( \frac{|x|^2}{ab} \Big)^{-2}\frac{1}{\log |x|}\le \frac{ab^2}{|x|^3}
$$
for  $|x| > 2L$, and by \eqref{f1}
\begin{equation}  \label{k0bis}
  |\mathcal R_1A(x)|\leq   \rho( R )^{-1} \int_{\nep^{-\alpha}}^{\nep^{\alpha}}\int_{|y|<L/2}    \frac{ab^2}{|x|^3}\,b^{-2}\frac{\di y\di b}{b}\lesssim \frac{a}{|x|^3},\, \qquad |x| > 2L.
\end{equation}
We proceed mainly as before, and define
\begin{equation*}
  S_0 = \Bbb R^2 \times [e^{-1-3\alpha/2}, e^{1+3\alpha/2} ].
\end{equation*}
The argument is split into two subcases, depending on the size of
$L\in [\nep^{2\alpha}, \nep^{8\alpha})$.

\noindent \textbf{Subcase (i):} $L\ge \nep^{3\alpha}$.
 Here we let
\begin{equation*}
  B_{0,j}=\{x\in S_0:\:2\nep^{j-1}L\le |x|\le 2\nep^{j}L\},  \qquad j = 1,2, \dots,
\end{equation*}
and  
\begin{equation*}
  B_{0,0}=\{x\in S_0:\: |x|\le 2L  \}.
\end{equation*}
Further, we define      
\begin{equation*}
 Z_{0,j}= [-2\nep^{j+2} L,2\nep^{j+2} L  ]^2 \times
[\nep^{-1-3\alpha/2-j/2},\, \nep^{1+3\alpha/2+j/2}] ,  \qquad j = 0,1, \dots.
\end{equation*}
Then $B_{0,j} \subset Z_{0,j} $, and since  $L\ge \nep^{3\alpha}$ the $Z_{0,j}$
 are   Calder\'on--Zygmund sets, centered at $(0,0,1)$ and of parameter $1+{3\alpha/2+j}/{2}$.
 The measures of these sets are given by
\begin{equation} \label{0meas}
  \rho(B_{0,j}) \sim \nep^{2j}L^2 \alpha  \qquad \mathrm{and}  \qquad
\rho(Z_{0,j}) \sim \nep^{2j} L^2 (\alpha + j).
\end{equation}
Now \eqref{k0bis} implies that for $x\in  B_{0,j},\; j\ge 1$
\begin{equation} \label{r1a0j}
|\mathcal R_1A(x)|\lesssim \frac{\nep^{3\alpha/2}}{\nep^{3j}L^3}
\end{equation}
As before, we define for
$j = 0,1, \dots$
\begin{equation}\label{m0j}
  m_{0,j} = \sum_{\ell = j}^\infty \int_{B_{0,\ell}} \mathcal R_1 A \,\di\rho.
\end{equation}
Since $\mathcal R_1 A \in L^1$,
 \eqref{vanint} implies   $\int_{S_0} \mathcal R_1 A \,d\rho = 0$,
so that $m_{0,0} = 0$.
Further,  \eqref{0meas} and \eqref{r1a0j} imply
\begin{equation*}
   |m_{0,j}| \lesssim  \frac{\nep^{3\alpha/2}\,\alpha}{\nep^{j}L},
\qquad j \ge 1.
\end{equation*}
We again let   $f_{0,j}=\rho(Z_{0,j})^{-1}\chi_{Z_{0,j}}$
 for $j = 0, 1,\dots$ and define
\begin{equation}\label{A0j}
  A_{0,j}=\mathcal R_1A\, \chi_{B_{0,j}} -m_{0,j}\, f_{0,j} + m_{0,j+1} \,f_{0,j+1}.
\end{equation}
Then
\begin{equation*}
 \sum_{j=0}^\infty A_{0,j}= \mathcal R_1A\, \chi_{\Omega_2} - m_{0,0}\, f_{0,0}
  =  \mathcal R_1A\, \chi_{\Omega_2}.
\end{equation*}

 Each  $A_{0,j}$ is supported in $Z_{0,j+1}$ and has integral 0.
Because of  \eqref{r1a0j}, $A_{0,j}$ is bounded for  $j\ge 1$ and is thus
a multiple of  a $(1,\infty)$-atom, with
\begin{equation*}
  \|A_{0,j}\|_{H^1}
 \lesssim \nep^{2j}L^2(\alpha+j)\,
\frac{\nep^{3\alpha/2}}{\nep^{3j} L^3}
+ \frac{\nep^{3\alpha/2}\,\alpha}{\nep^{j}L}
\lesssim  \frac{ j}{\nep^{j}},\,  \qquad j\ge 1,
\end{equation*}
since $L \ge \nep^{3\alpha}$ in this subcase. Thus
\begin{equation*}
\| \sum_1^\infty  A_{0,j}\|_{H^1} \lesssim 1.
\end{equation*}
It only remains to
consider $A_{0,0}$, which  need not be bounded. We use now the boundedness of $\mathcal R_1$
on $L^2$, which implies
\begin{equation}\label{A00}
\begin{aligned}
   \|A_{0,0}\|_{2} &\lesssim
 \| \mathcal R_1A\|_{2}
+   |m_{0,1}|\, \rho(Z_{0,1})^{-1/2}\\
& \lesssim
 \|  A\|_{2} + \frac{\nep^{3\alpha/2}{\alpha}}{L}\, \rho(Z_{0,1})^{-1/2}
\lesssim  \rho(Z_{0,1})^{-1/2}.
\end{aligned}
\end{equation}
Thus $A_{0,0}$ is a multiple of a  $(1,2)$-atom.

Summing up,  we conclude that
$\mathcal R_1A\, \chi_{\Omega_2}\in {H^1}$ with bounded norm.
                        This ends Subcase (i).

\medskip

\noindent \textbf{Subcase (ii):} $L < \nep^{3\alpha}$.  Here we divide
$S_0$ into two slices, defining
\begin{equation*}
  S_{0-} = \Bbb R^2 \times [e^{-1-3\alpha/2},\, e^{\alpha/2} ]
\end{equation*}
and
\begin{equation*}
  S_{0+} = \Bbb R^2 \times [e^{\alpha/2},\, e^{1+3\alpha/2} ].
\end{equation*}
Further, we let for $j = 1,2, \dots$
\begin{equation*}
  B_{0-,j}=\{x\in S_{0-}:\:2\nep^{j-1}L\le |x|\le 2\nep^{j}L\}
\end{equation*}
and
\begin{equation*}
 B_{0+,j}=\{x\in S_{0+}:\:2\nep^{j-1}L\le |x|\le 2\nep^{j}L\},
\end{equation*}
but for $j=0$, as before,
\begin{equation*}
 B_{0\pm,0}=\{x\in S_{0\pm}:\: |x|\le 2L\}.
\end{equation*}
Then
\begin{equation*}
 Z_{0-,j}:= [-2\nep^{j+2} L,2\nep^{j+2} L  ]^2 \times
[\nep^{-1-3\alpha/2-j/2},\, \nep^{\alpha/2+j/2}] ,  \qquad j = 0,1, \dots,
\end{equation*}
are Calder\'on--Zygmund sets of center $(0,0,\nep^{-(\alpha+1)/2})$ and parameter $\alpha+(j+1)/{2}$ containing $B_{0-,j}$. Almost similarly,
\begin{equation*}
 Z_{0+,j}:= [-2\nep^{j+2} L,2\nep^{j+2} L  ]^2 \times
[\nep^{\alpha/2-j/2}, \,\nep^{1+3\alpha/2+j/2}] ,  \qquad j = 0,1, \dots,
\end{equation*}
are Calder\'on--Zygmund sets of center $(0,0,\nep^{\alpha+1/2})$ and parameter $(\alpha+1+j)/{2}$ containing $B_{0+,j}$. The measures of these sets are
\begin{equation*}
  \rho(B_{0\pm,j}) \sim \nep^{2j}L^2 \alpha  \qquad \mathrm{and}  \qquad
\rho(Z_{0\pm,j}) \sim \nep^{2j} L^2 (\alpha + j).
\end{equation*}
 For $x\in  B_{0-,j}$ with $j\ge 1$, the estimate \eqref{k0bis} implies
\begin{equation}\label{RA-}
|\mathcal R_1A(x)|\lesssim \frac{\nep^{\alpha/2}}{\nep^{3j}L^3}
\end{equation}
and for $x\in  B_{0+,j}$ with $j\ge 1$
\begin{equation}\label{RA+}
|\mathcal R_1A(x)|\lesssim \frac{\nep^{3\alpha/2}}{\nep^{3j}L^3}.
\end{equation}
Proceeding as before,  for $j = 0,1, \dots$ we define
 $f_{0\pm,j}=\rho(Z_{0\pm,j})^{-1}\chi_{Z_{0\pm,j}}$ and
\begin{equation*}
  m_{0\pm,j} = \sum_{\ell = j}^\infty \int_{B_{0\pm,\ell}} \mathcal R_1 A \,\di\rho.
\end{equation*}
 and for $j\ge 1$ also
\begin{equation}\label{A0j+-}
  A_{0\pm,j}=\mathcal R_1A\, \chi_{B_{0\pm,j}} -m_{0\pm,j}\, f_{0\pm,j}
+ m_{0\pm,j+1}\, f_{0\pm,j+1}.
\end{equation}
These $A_{0\pm,j}$ are multiples of $(1,\infty)$-atoms, because they are bounded and
supported in $Z_{0\pm,j+1}$ and have vanishing integrals. Since
\eqref{RA-} and \eqref{RA+} imply
\begin{equation*}
   |m_{0\pm,j}| \lesssim  \frac{\nep^{3\alpha/2}\alpha}{\nep^{j}L},
\qquad j\ge 1,
\end{equation*}
one finds that
\begin{equation*}
  \|A_{0\pm,j}\|_{H^1}
 \lesssim \frac{\nep^{3\alpha/2} (\alpha+j)}{\nep^{j}L}
+ \frac{\nep^{3\alpha/2}\alpha}{\nep^{j}L} \lesssim  \frac{ j}{\nep^{j}},\,  \qquad j\ge 1.
\end{equation*}
The sum of these $A_{0\pm,j}$ is thus an $H^1$ function, and
\begin{equation*}
  \sum_{j=1}^\infty A_{0\pm,j} =  \mathcal R_1A\,\chi_{\cup_1^\infty B_{0\pm,j}}
- \int_{\cup_1^\infty B_{0\pm,j}}  \mathcal R_1A\,\di\rho.
\end{equation*}
These two  sums form a large part of the desired atomic decomposition of
$\mathcal R_1A\chi_{S_0}$. What is missing is
\begin{multline}\label{splusminus}
  \mathcal R_1A\chi_{S_0} -  \sum_{j=1}^\infty A_{0-,j}-  \sum_{j=1}^\infty A_{0+,j} \\
 = \mathcal R_1A\,\chi_{B_{0-,0}} + \mathcal R_1A\,\chi_{B_{0+,0}}
+  \int_{\cup_1^\infty B_{0-,j}} \mathcal R_1A\,\di\rho + \int_{\cup_1^\infty B_{0+,j}} \mathcal R_1A\,\di\rho.
\end{multline}
Comparing with Subcase (i), we have   $S_{0-}\cup S_{0+} = S_{0} $
and  $B_{0-,j}\cup B_{0+.j} = B_{0,j} $ for $j = 0,1,\dots$,
and we observed there that
  $\int_{S_{0}} \mathcal R_1A\,\di\rho = 0$.
 It follows that the last two terms in the right-hand side of \eqref{splusminus}
sum up to $-\int_{B_{0,0}} \mathcal R_1A\,\di\rho$, and
the first two terms amount to $ \mathcal R_1A\,\chi_{B_{0,0}}$.
 This means that the expression
in  \eqref{splusminus} coincides with the atom multiple $A_{0,0}$ from
 Subcase (i), and the estimate \eqref{A00} holds in both subcases. We thus have an atomic decomposition of
$\mathcal R_1A\,\chi_{S_0}$ in  Subcase (ii), which ends the proof of
Case I in Theorem \ref{t: H1bdd}.

 \bigskip

\noindent{\bf{Case II: intermediate atom.}} In this case $1\leq \alpha<20$ and $\nep^{2\alpha}\leq L<\nep^{8\alpha}$. We consider $A$ as a function supported in the larger Calder\'on--Zygmund set $[-\tilde L/2,\tilde L/2]^2\times [\nep^{-20},\nep^{20}]$, where $\tilde L=\max\{L,\nep^{40}\}$. Considered with this support, $A$ will be a multiple of a large atom, and the preceding argument applies.

\bigskip

\noindent {\bf{Case III: small atom.}} Here $\alpha <1$ and $\nep^2\alpha\leq L<\nep^8\alpha$. We use a slightly different decomposition of the space into regions, as follows
$$
\begin{aligned}
\tilde \Omega_1&=\{(x_1,x_2,a):\,a\leq \nep^{-2} \}\,,\\
\tilde \Omega_2&=\{(x_1,x_2,a):\,  \nep^{-2}<a<\nep^2\}\,,\\
\tilde \Omega_3&=\{(x_1,x_2,a):\,a\geq \nep^{2} \}\,.
\end{aligned}
$$
Starting with $\tilde \Omega_1$, we let  $x\in \tilde\Omega_1$ and $y = (y_1,y_2,b)\in R$.
Then $b^{-1}a < e^{-1}$, and
\eqref{d-piccolo}  shows that
$$
|k_1(y^{-1}x)|\lesssim   \frac{a}{1+|x|^3} \qquad {\rm{and}} \qquad |\mathcal R_1A(x)|\lesssim  \frac{a}{1+|x|^3} \,.
$$
We can construct an atomic decomposition of
$\mathcal R_1A \chi_{\tilde \Omega_1}$ using the following sets
for $k<0$ and $j\geq 0$:
$$
\begin{aligned}
B_{k,j}&= \{ 2\nep^{j-1} \leq |x|\leq 2\nep^j \}\times   [\nep^{-2^{|k|}-2 } , \nep^{-2^{|k|-1}-2 }]\,,\qquad j\geq 1\,,\\    
B_{k,0}& =\{ |x|\leq 2 \} \times [\nep^{-2^{|k|}-2} , \nep^{-2^{|k|-1}-2}]\,,\\
Z_{k,j}&=[-2 \nep^j ,2\nep^j ]^2\times [\nep^{-3\cdot 2^{|k|-1}-2 -j/2} , \nep^{-2^{|k|-1}-2 +j/2}],\qquad j\geq 0\,.
\end{aligned}
$$
The $Z_{k,j}$ are Calder\'on--Zygmund sets, centered at $(0,0,\nep^{-2\cdot 2^{|k|-1}-2})$ and with parameters $2^{|k|-1}+j/{2}$. The relevant  measures are $\rho(B_{k,j})\sim \nep^{2j}2^{|k|}$         and $\rho(Z_{k,j})\sim \nep^{2j}(2^{|k|}+j)$. We then have in the set $B_{k,j}$

$$
|\mathcal R_1 A|\lesssim \frac{\nep^{-2^{|k|-1}}}{\nep^{3j}},\,
\qquad k<0, \quad j \ge 0.
$$
From here, we proceed as in the case of the region $\Omega_1$ for a large atom. The details are left to the reader, since the construction is quite similar. We will have $\|\mathcal R_1A\chi_{\tilde \Omega_1}\|_{H^1}\lesssim 1$\,.

In the case of $\tilde \Omega_3$, we argue as we did to obtain
the estimates \eqref{R_1A'} and \eqref{R_1A''}, to show that
for $x=(x_1,x_2,a)\in \tilde \Omega_3$
$$
|\mathcal R_1A(x)|\lesssim \frac{La}{(a^2+|x|^2)^2 \log a  }+\frac{a|x_1|\alpha}{ (a^2+|x|^2)^2   ( \log  a)^2  }, \qquad |x|> 2,
$$
and \eqref{cgrande} leads to
$$
|\mathcal R_1A(x)|\lesssim \frac{1}{a^3\log a}, \qquad |x| \le 2.
$$
We now construct an atomic decomposition of $\mathcal R_1A \chi_{\tilde \Omega_3}$, using the following sets for $k\geq 1$ and $j\geq 0$:
$$
\begin{aligned}
B_{k,j}&= \{ \nep^{k+j+1} \leq |x|\leq \nep^{k+j+2} \}\times  [\nep^{k+1} , \nep^{k+2}] \,,\quad j\geq 1\,,\\
B_{k,0}&=\{|x|\leq e^{k+2}\}\times [\nep^{k+1} , \nep^{k+2}]\,,\\
Z_{k,j}&=   \{   |x|\leq \nep^{k+j+2} \}   \times  [\nep^{k+1-j/2} , \nep^{k+2+j/2}] \,.
\end{aligned}
$$
Here each $Z_{k,j}$ is a Calder\'on--Zygmund set of center $(0,0,\nep^{k+3/2})$ and parameter $(j+1)/{2}$ containing $B_{k,j}$.
The measures satisfy  $\rho(B_{k,j})\sim \nep^{2k+2j}$ and
 $\rho(Z_{k,j})\sim \nep^{2k+2j}(j+1)$. In $B_{k,j}$ we will have
$$
|\mathcal R_1 A|\lesssim \frac{1}{\nep^{3k+4j}k} +   \frac{1}{\nep^{2k+3j}k^2}\,.
$$

Again, we can now follow the procedure used for a large atom in the region
 $\Omega_1$.
In particular, $f_k$, \hskip2pt $m_{k,j}$ and $A_{k,j}$ will be  as there,
although $k$ is now positive and the  $B_{k,j}$ and  $Z_{k,j}$ are those
we just defined. The relevant estimates will now be
\begin{align}
  |m_{k,j}| \lesssim\, & \sum_{\ell = j}^\infty \left( \frac{1}{\nep^{3k+4\ell}k} \,
\nep^{2k+2\ell} + \frac{1}{\nep^{2k+3\ell}k^2} \,\nep^{2k+2\ell} \right)
\notag \\ \notag
\lesssim &\, \nep^{-k-2j} + \nep^{-j} \frac{1}{k^2}
\end{align}
and
\begin{align}
 & \| A_{k,j} \|_{H^1} \notag \\ &\lesssim \, \frac{1}{\nep^{3k+4j}\,k}\, \nep^{2k+2j}\,(j+1)
+   \frac{1}{\nep^{2k+3j}\,k^2} \,\nep^{2k+2j}\,(j+1)
+ \nep^{-k-2j} + \nep^{-j} \frac{1}{k^2}
 \notag \\
\notag & \lesssim  \:\nep^{-k-2j}\,(j+1) + \nep^{-j}\, \frac{j+1}{k^2}.
\end{align}
These quantities  are summable over $j \ge 0$ and $k\ge 1$,
and so
 $\|\mathcal R_1A\chi_{\tilde \Omega_3 }\|_{H^1}\lesssim 1$.

\smallskip

To treat the region $\tilde{\Omega}_2$ we will construct a sequence of sets $B_{0,j}$ expanding from $R$ in all three coordinates until they reach essentially unit size. Then they will expand only in the coordinates $x_1,x_2$. More precisely, we never let the $a$ width of any $B_{0,j}$ be larger than
$[\nep^{-2},   \nep^{2}]$, so that $B_{0,j}$ stays in $\tilde{\Omega}_2$. A sequence of Calder\'on--Zygmund sets $Z_{0,j}$ will be defined accordingly. It is not restrictive to suppose that $\alpha=2^{j_0}$, for some integer $j_0 \le 0$. The definition of a Calder\'on--Zygmund set
then implies that $e^2 \le 2^{-j_0} L < e^8$. The  $B_{0,j}$ and the  $Z_{0,j}$  will be defined for
$j = j_0+1, j_0+2,\dots$, as follows.

We start by setting $B_{0,j_0} = \emptyset$
  and recursively    for $ j \ge j_0+1 $
\begin{equation*}
  B_{0,j}=\left(\{|x|< 2^{j-j_0-1}L\}\times [\nep^{-\min(2,2^{j})},   \nep^{\min(2,2^{j})}]\right)\setminus  B_{0,j-1}\,,
\end{equation*}

For  $  j_0+1 \le j \le 0$ (which occurs only if $j_0 < 0 $ ), we let
$$
Z_{0,j}=[-2^{j-j_0-1}L,2^{j-j_0-1}L]^2\times [\nep^{-2^{j}},   \nep^{2^{j}}]\,,
$$
but for $j > 0$
$$
Z_{0,j}=[-2^{j-j_0+2}L,2^{j-j_0+2}L]^2\times [\nep^{-2-j/2},   \nep^{2+j/2}]\,.
$$
Then $B_{0,j} \subset Z_{0,j}$ for each $j \ge j_0+1$, and the $Z_{0,j}$ are Calder\'on--Zygmund sets
centered at $(0,0,1)$ and of parameter $\min(2,2^{j})+j_+/2$. The measures of these sets are
$\rho(B_{0,j}) \sim 2^{2j+j_-}$ and $\rho(Z_{0,j}) \sim 2^{2j}(2^{j_-}+j_+)$.
Here $j_+ = \max(j,0)$ and $j_- = \min(j,0)$.

\vskip2pt    

Suppose now that $j_0+2 \le j\leq 0$. We bound $\mathcal R_1A$ in the set $B_{0,j}$ by means of  estimates similar to those for the Riesz transforms in the Euclidean setting. Observe first that each point $x=(x_1,x_2,a)\in B_{0,j}$ is at some distance from any point $y$ in $R$, so that $r(y^{-1}x)\sim r(x)$. Then simple computations together with formulas \eqref{k1},  \eqref{derk1}, \eqref{k0}  and  \eqref{derivativedistance} show that
$$
\begin{aligned}
|k_1(y^{-1}x)|&\,\lesssim\, r(x)^{-3}\,,\\   
 |X_ik_1(y^{-1}x)|&\,\lesssim \,r(x)^{-4}\,,\qquad  i=0,1,2\,.
\end{aligned}
$$
Notice that when $r(x)$ is small, it is essentially the Euclidean distance from $x$ to $e$; indeed
 $r(x) \sim \sqrt{x_1^2+x_2^2+(\log a)^2}$ as seen from \eqref{metrica}.
From the Mean Value Theorem and the fact that $\delta(y) = 1 + \mathcal O(\alpha)$ here,
it now follows that for  $x\in B_{0,j}$
\begin{align*}                 
|\mathcal R_1A(x)|&\lesssim \int_R |A(y)||k_1(y^{-1}x)\delta(y)-k_1(x)|\dir(y) \\
&\lesssim  \int_R |A(y)||k_1(y^{-1}x)-k_1(x)|\dir(y)
+ \alpha  \int_R |A(y)||k_1(y^{-1}x) \dir(y)\\
& \lesssim \frac{\alpha}{r(x)^{4}} + \frac{\alpha}{r(x)^{3}}\lesssim \alpha 2^{-3j}\,.
\end{align*}

Suppose now instead that  $j>0$. Then \eqref{k0bis}  shows that for each $x\in B_{0,j}$
$$
|\mathcal R_1A(x)|\lesssim \frac 1{|x|^3}\lesssim 2^{-3j}\,.
$$
But we can also apply the Mean Value Theorem as in the preceding case.   Then we need \eqref{stimaderk1}
 to estimate $X_i k_1$ for $i = 1,2$, and by means of  \eqref{k1} and  \eqref{derivativedistance} one
can verify that the same estimate holds also for $X_0 k_1$. The result will be
$$
|\mathcal R_1A(x)|\lesssim \frac \alpha{|x|^2\log|x|}\lesssim \frac{\alpha 2^{-2j}}j\,.
$$

We now argue essentially as in Subcase (i) above. In particular, we define $m_{0,j}$ by \eqref{m0j} and $A_{0,j}$ by \eqref{A0j}, for all $j \ge j_0 + 1$. Then
$\sum_{j=j_0+1}^{\infty}A_{0,j}=\mathcal R_1\,\chi_{\tilde \Omega_2}$, and
$m_{0,j_0+1} =\int_{B_{\tilde \Omega_2}} \mathcal R_1 A\,d\rho = 0$.

Using these estimates for $\mathcal R_1A$, we  estimate $m_{0,j}$, first when $ j_0 + 1 < j \le 0$. Then
\begin{align*}
  |m_{0,j}|\, \lesssim & \,\sum_{\ell=j}^\infty \int_{B_{0,j}} \mathcal R_1 A\,d\rho \lesssim
  \sum_{\ell=j}^0   \alpha 2^{-3\ell}\, 2^{3\ell}  +  \sum_{\ell=1}^\infty
 \min\left( 2^{-3\ell}\, 2^{2\ell} \ell, \frac{\alpha 2^{-2\ell}}\ell \,2^{2\ell} \ell\right)\\
  \lesssim & \,\alpha\, (1+ \log1/\alpha),
\end{align*}
the last step since the number of terms in the finite sum is at most
$|j_0| \lesssim 1+\log1/\alpha$ and the last sum is easy to control.

If $j>0$ we have instead
\begin{equation*}
  |m_{0,j}| \lesssim   \sum_{\ell=j}^\infty 2^{-3\ell}\, 2^{2\ell} \ell \lesssim  2^{-j} j .
\end{equation*}


    Our estimates for $\mathcal R_1A$  and $\rho(Z_{0,j})$  now show that
$$
\begin{aligned}
\|A_{0,j}\|_{H^1}&\lesssim \alpha \,(1+ \log1/\alpha)\,, \qquad j_0+1<j\leq 0\,,\\
\|A_{0,j}\|_{H^1}&\lesssim 2^{-j}j\qquad j>0\,.
\end{aligned}
$$
Finally, we use the $L^2$-boundedness of $\mathcal R_1$ as in \eqref{A00} to deduce that $A_{0,j_0}$ is a multiple of a $(1,2)$-atom. We can now sum the $H^1$-norms of all the $A_{0,j}$ and obtain $\|\mathcal R_1A\,\chi_{\tilde \Omega_2}\|_{H^1}\lesssim 1$.

\smallskip

This concludes the proof of Theorem \ref{t: H1bdd}.



\begin{thebibliography}{CDKR1}





\bibitem{ADY} J.~Ph.~Anker, E.~Damek, C.~Yacoub,
                  \emph{Spherical Analysis on harmonic $AN$ groups},
                  Ann. Scuola Nom. Sup. Pisa Cl. Sci. 23 (1996), 643--679.







\bibitem{AMR} P. Auscher, A. McIntosh, E. Russ, \emph{Hardy spaces of differential forms and Riesz transforms on Riemannian manifolds}, C. R. Math. Acad. Sci. Paris 344 (2007), no. 2, 103--108.

\bibitem{BDG} J. Betancor, J. Dziuba\'nski, G. Garrig\'os, \emph{Riesz transform characterization of Hardy spaces associated with certain Laguerre expansions}, Tohoku Math. J. (2) 62 (2010), no. 2, 215--231.

\bibitem{CaMaMe} A. Carbonaro, G. Mauceri, S. Meda, \emph{$H^1$ and $BMO$ for certain locally doubling metric measure spaces}, Ann. Sc. Norm. Super. Pisa Cl. Sci. (5) 8 (2009), no. 3, 543--582.

\bibitem{CMM} A. Carbonaro, A. McIntosh, A.J. Morris, \emph{Local Hardy spaces of differential forms on Riemannian manifolds}, J. Geom. Anal. 23 (2013), no. 1, 106--169.




\bibitem{CeMe} D. Celotto, S. Meda, \emph{On the analogue of the Fefferman-Stein theorem on graphs with the Cheeger property}, Ann. Mat. Pura Appl. (4) 197 (2018), no. 5, 1637--1677.

\bibitem{CW1} R. R. Coifman, G. Weiss,
                 \emph{Analyse harmonique non commutative sur certains espaces homogenes},
                 Lecture Notes in Mathematics 242, Springer (1971).

\bibitem{CW2} R. R. Coifman, G. Weiss,
                   \emph{Extensions of Hardy spaces and their use in Analysis},
                   Bull. Am. Math. Soc. 83 (1977), 569--645.


%
%

\bibitem{CGHM} M.~Cowling, S.~Giulini, A.~Hulanicki, G.~Mauceri,  {\emph{Spectral multipliers for a distinguished laplacian on   certain groups of exponential growth}}, Studia Math. {111} (1994), 103--121.


















\bibitem{DP} J. Dziuba\'nski, M. Preisner, \emph{Riesz transform characterization of Hardy spaces associated with Schr\"odinger operators with compactly supported potentials}, Ark. Mat. 48 (2010), no. 2, 301--310.

\bibitem{D} J. Dziuba\'nski, \emph{Riesz transforms characterizations of Hardy spaces $H^1$ for the rational Dunkl setting and multidimensional Bessel operators}, J. Geom. Anal. 26 (2016), no. 4, 2639--2663.

\bibitem{GS} S.~Giulini, P.~Sj\"ogren,
                 \emph{A note on maximal functions on a solvable Lie group},
                 Arch. Math. (Basel) 55 (1990), 156--160.

%
%

\bibitem{HS} W.~Hebisch, T.~Steger,
                 \emph{Multipliers and singular integrals on expo\-nen\-tial growth groups},
                 Math. Z. 245 (2003), 37--61.


\bibitem{HLMMY} S. Hofmann, G. Lu, D. Mitrea, M. Mitrea, L. Yan,
\emph{Hardy spaces associated to non-negative self-adjoint operators satisfying Davies-Gaffney estimates},
Mem. Amer. Math. Soc. 214 (2011), no. 1007, vi+78 pp.

\bibitem{LMSTV} M. Levi, A. Martini, F. Santagati, A. Tabacco, M. Vallarino,
     \emph{Riesz transform for a flow {L}aplacian on homogeneous trees}, {J. Fourier Anal. Appl.} {29} (2023), no. 2, Paper No. 15, 29 pp.



\bibitem{LoVa} N. Lohou\'e, N. Varopoulos, \emph{Remarques sur les transform\`ees de Riesz sur les groupes de Lie nilpotents}, C. R. Acad. Sci. Paris S\'er. I Math. 301 (1985), no. 11, 559--560.




\bibitem{MaRu} M. Marias, E. Russ, \emph{$H^1$-boundedness of Riesz transforms and imaginary powers of the Laplacian on Riemannian manifolds},
Ark. Mat. 41 (2003), no. 1, 115--132.









\bibitem{Ma} A. Martini, \emph{Riesz transforms on  ax+b  groups}, J. Geom. Anal. 33 (2023), no. 7, 53 pp.


\bibitem{MaMeV} A. Martini, S. Meda, M. Vallarino, \emph{A family of Hardy type spaces on nondoubling manifolds}, Ann. Mat. Pura Appl. 199 (2020), no. 5, 2061--2085.


\bibitem{MaMeVaVe} A. Martini, S.  Meda, M. Vallarino, G. Veronelli, \emph{Inclusions and noninclusions of Hardy type spaces on certain nondoubling manifolds},
J. Funct. Anal. 286 (2024), no. 3, 110240.

\bibitem{MSV} A. Martini, F. Santagati, M. Vallarino, \emph{Heat kernel and Riesz transform for the flow Laplacian on homogeneous trees}, to appear in N. Arcozzi, M. M. Peloso and A.
Tabacco, (New Trends in) Complex and Fourier Analysis, Springer INdAM Series, Springer. arXiv:2210.07148


\bibitem{MSTV} A. Martini, F. Santagati, A. Tabacco, M. Vallarino, \emph{Riesz transform and spectral multipliers for the flow Laplacian on nonhomogeneous trees}, arXiv:2310.09113


\bibitem{MaVa} A. Martini. M. Vallarino, \emph{Riesz transforms on solvable extensions of stratified groups}, Studia Math. 259 (2021), no. 2, 175--200.



\bibitem{MMV} G. Mauceri, S. Meda, M. Vallarino, \emph{Higher order Riesz transforms on noncompact symmetric spaces}, J. Lie Theory 28 (2018), no. 2, 479--497.

 \bibitem{MMS} G. Mauceri, S. Meda, P. Sj\"ogren, \emph{Endpoint estimates for first-order Riesz transforms associated to the Ornstein-Uhlenbeck operator}, Rev. Mat. Iberoam. 28 (2012), no. 1, 77--91.


 \bibitem{MeVe} S. Meda, G. Veronelli, \emph{Local Riesz transform and local Hardy spaces on Riemannian manifolds with bounded geometry},
J. Geom. Anal. 32 (2022), no. 2, Paper No. 55, 57 pp.

\bibitem{Ru} E. Russ, \emph{$H^1-L^1$ boundedness of Riesz transforms on Riemannian manifolds and on graphs}, Potential Anal. 14 (2001), no. 3, 301--330.

\bibitem{SC} L. Saloff-Coste, \emph{Analyse sur les groupes de Lie \`a croissance polynomiale}, Ark. Mat. 28 (1990), 315--331.

\bibitem{S} F. Santagati, \emph{Hardy spaces on homogeneous trees with flow measures}, {J. Math. Anal. Appl.} {510} (2022), 23 pp.


\bibitem{SV1} P.~Sj\"ogren, M. Vallarino,
               \emph{ Boundedness from $H^1$ to $L^1$ of Riesz transforms on a Lie group of exponential growth},  Ann. Inst. Fourier
58 (2008), 1117--1151.

\bibitem{SV2} P.~Sj\"ogren, M. Vallarino,
               \emph{Heat maximal function on a Lie group of exponential growth},  Ann. Acad. Sci. Fenn. Math. 37 (2012), no. 2, 491--507.

\bibitem{St} E.M.~Stein,
               \emph{Harmonic Analysis},
               Princeton University Press (1993).



\bibitem{V} M.~Vallarino,
             \emph{Spaces $H^1$ and $BMO$ on exponential growth groups}, Collect. Math. 60 (2009), 277--295.

\bibitem{V2} M.~Vallarino, \emph{Atomic and maximal Hardy spaces on a Lie group of exponential growth}, Trends in harmonic analysis, 409--424, Springer INdAM Ser., 3, Springer, Milan, 2013.



\end{thebibliography}
\end{document}